\documentclass[twoside,11pt]{article}

\usepackage{blindtext}

\usepackage{amsthm}
\usepackage[abbrvbib, preprint]{jmlr2e}
\usepackage{amsmath,amssymb,amssymb}
\usepackage[capitalize]{cleveref}
\usepackage{comment}
\usepackage{tikz}
\usetikzlibrary{calc}
\usetikzlibrary{positioning,arrows.meta,shapes,arrows}
\usepackage{textcomp}
\tikzset{%
	block/.style    = {draw, thick, rectangle, minimum height = 3em,
			minimum width = 3em},
	sum/.style      = {draw, circle, node distance = 2cm}, 
	input/.style    = {coordinate}, 
	output/.style   = {coordinate} 
}
\newcommand{\suma}{\Large$+$}
\newcommand{\inta}{\Large$\int$}
\newcommand{\muma}{\Large$\times$}

\usepackage{amsfonts}%
\usepackage{amssymb}%
\usepackage{mathrsfs}
\usepackage{overpic}
\usepackage{mathtools}
\usepackage{enumitem}
\usepackage{bm}
\usepackage{multirow}
\usepackage{algorithm,float}
\usepackage{algpseudocode}

\newcommand{\R}{\mathbb{R}}

\newcommand{\Lop}{\mathcal{L}}
\newcommand{\diam}{\textup{diam}}
\newcommand{\dist}{\textup{dist}}

\renewcommand{\d}{\,\textup{d}}

\DeclareMathOperator{\vdiv}{div}

\newcommand{\A}{\mathcal{A}}
\DeclareMathOperator{\U}{\mathcal{U}}
\DeclareMathOperator{\V}{\mathcal{V}}
\renewcommand{\div}{\mathrm{div}}

\usepackage{lastpage}
\jmlrheading{XX}{2023}{1-\pageref{LastPage}}{; Revised }{}{21-0000}{Boull\'e and Townsend}

\ShortHeadings{A Mathematical Guide to Operator Learning}{Boull\'e and Townsend}
\firstpageno{1}

\usepackage{xcolor}
\definecolor{myred}{RGB}{238,33,13} 
\definecolor{mygreen}{RGB}{30,176,1} 
\definecolor{myblue}{RGB}{1, 161, 255} 
\definecolor{myyellow}{RGB}{255, 240, 87} 
\definecolor{myorange}{RGB}{242, 115, 0} 

\begin{document}

\title{A Mathematical Guide to Operator Learning}

\author{\name Nicolas Boull\'{e} \email nb690@cam.ac.uk \\
\addr Department of Applied Mathematics and Theoretical Physics\\
University of Cambridge\\
Cambridge, CB3 0WA, UK
\AND
\name Alex Townsend \email townsend@cornell.edu \\
\addr Department of Mathematics\\
Cornell University\\
Ithaca, NY 14853, USA
}

\editor{}

\maketitle
\vskip 0.2in

\begin{abstract}
    Operator learning aims to discover properties of an underlying dynamical system or partial differential equation (PDE) from data. Here, we present a step-by-step guide to operator learning. We explain the types of problems and PDEs amenable to operator learning, discuss various neural network architectures, and explain how to employ numerical PDE solvers effectively. We also give advice on how to create and manage training data and conduct optimization. We offer intuition behind the various neural network architectures employed in operator learning by motivating them from the point-of-view of numerical linear algebra.
\end{abstract}

\vskip 0.2in


\tableofcontents

\section{Introduction}

The recent successes of deep learning~\citep{lecun2015deep} in computer vision~\citep{krizhevsky2012imagenet}, language model~\citep{brown2020language}, and biology~\citep{jumper2021highly} have caused a surge of interest in applying these techniques to scientific problems. The field of scientific machine learning (SciML)~\citep{karniadakis2021physics}, which combines the approximation power of machine learning (ML) methodologies and observational data with traditional modeling techniques based on partial differential equations (PDEs), sets out to use ML tools for accelerating scientific discovery.

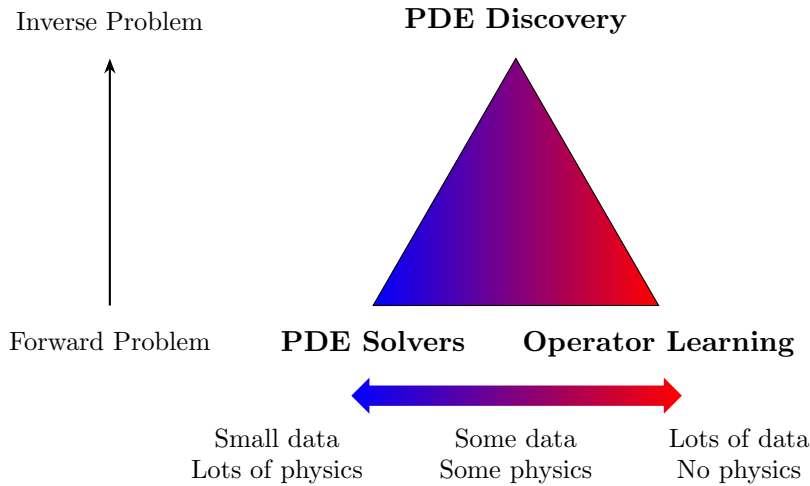
\begin{figure}[htbp]
    \centering
    \makeatletter
\newcommand{\gettikzxy}[3]{%
    \tikz@scan@one@point\pgfutil@firstofone#1\relax
    \edef#2{\the\pgf@x}%
    \edef#3{\the\pgf@y}%
}
\makeatother
\usetikzlibrary{shapes.arrows}
\usetikzlibrary{calc}

\begin{tikzpicture}[nodes={align=center}]
    \path node[regular polygon, regular polygon sides=3,
        inner sep=2em, fill opacity=0] (3gon) {};

    \node[above, font=\bfseries] (x1) at ([yshift=2mm]3gon.corner 1){{PDE Discovery}};

    \node[below, font=\bfseries] (x2) at ([yshift=-2mm]3gon.corner 2){{PDE Solvers}};

    \node[below, font=\bfseries] (x3) at ([yshift=-2mm]3gon.corner 3){{Operator Learning}};

    \draw[fill, top color=blue, bottom color=red, shading=axis, shading angle=90.02]
    (3gon.corner 1) -- (3gon.corner 2) -- (3gon.corner 3) -- cycle;
    \gettikzxy{(x2)}{\ax}{\ay}
    \gettikzxy{(x1)}{\bx}{\by}
    \gettikzxy{(x3)}{\ex}{\ey}

    \gettikzxy{(3gon.corner 2)}{\cx}{\cy}
    \gettikzxy{(3gon.corner 1)}{\dx}{\dy}
    \gettikzxy{(3gon.corner 3)}{\fx}{\fy}

    \node[font=\small] at (\ax-35mm,\ay){Forward Problem};
    \node[font=\small] at (\ax-35mm,\by){Inverse Problem};
    \draw[-{Stealth[length=2mm]}, thick] (\cx-35mm,\cy) -- (\cx-35mm, \dy);


    \node[below, font=\small] at ([yshift=-15mm]$(3gon.corner 2)!0.5!(3gon.corner 3)$){Some data\\Some physics};

    \node[below left, font=\small, align=center] at ([xshift=0mm,yshift=-15mm]3gon.corner 2) {Small data\\Lots of physics};

    \node[below right, font=\small, align=center] at ([xshift=0mm,,yshift=-15mm]3gon.corner 3) {Lots of data\\No physics};

    \node [double arrow,left color=blue,right color=red,
        double arrow head extend=3pt,transform shape,minimum height=\ex-\ax+6mm,anchor=west]
    at ([xshift=-3mm,yshift=-12mm]3gon.corner 2){};
\end{tikzpicture}
    \caption{Illustrating the role of operator learning in SciML. Operator learning aims to discover or approximate an unknown operator $\A$, which often corresponds to the solution operator of an unknown PDE. In contrast, PDE discovery aims to discover coefficients of the PDE itself, while PDE solvers aim to solve a known PDE using ML techniques.}
    \label{fig_schematic}
\end{figure}

SciML techniques can roughly be categorized into three main areas: (1) PDE solvers, (2) PDE discovery, and (3) operator learning (see~\cref{fig_schematic}). First, PDE solvers, such as physics-informed neural networks (PINNs)~\citep{raissi2019physics,lu2021deepxde,cuomo2022scientific,wang2023expert}, the deep Galerkin method~\citep{sirignano2018dgm}, and the deep Ritz method~\citep{yu2018deep}, consist of approximating the solution a known PDE by a neural network by minimizing the solution's residual. At the same time, PDE discovery aims to identify the coefficients of a PDE from data, such as the SINDy approach~\citep{brunton2016discovering,champion2019data}, which relies on sparsity-promoting algorithms to determine coefficients of dynamical systems. There are also symbolic regression techniques, such as AI Feynman introduced by~\citep{udrescu2020ai,udrescu2020ai2} and genetic algorithms~\citep{schmidt2009distilling,searson2010gptips}, that discover physics equations from experimental data.

Here, we focus on the third main area of SciML, called operator learning~\citep{lu2021learning,kovachki2023neural}. Operator learning aims to discover or approximate an unknown operator $\A$, which often takes the form of the solution operator associated with a differential equation. In mathematical terms, the problem can be defined as follows. Given pairs of data $(f,u)$, where $f\in \U$ and $u\in \V$ are from function spaces on a $d$-dimensional spatial domain $\Omega\subset\R^d$, and a (potentially nonlinear) operator $\A:\U\to\V$ such that $\A(f)=u$, the objective is to find an approximation of $\A$, denoted as $\hat{\A}$, such that for any new data $f'\in \U$, we have $\hat{\A}(f')\approx \A(f')$. In other words, the approximation should be accurate for both the training and unseen data, thus demonstrating good generalization.

This problem is typically approached by representing $\hat{\A}$ as a neural operator, which is a generalization of neural networks as the inputs and outputs are functions, not vectors. After discretizing the functions at sensor points $x_1,\ldots,x_m\in\Omega$, one then parametrizes the neural operator with a set of parameters $\theta\in \R^N$, which could represent the weights and biases of the underlying neural network. Then, one typically formulates an optimization problem to find the best parameters:
\begin{equation} \label{eq_loss_function}
    \min_{\theta\in \R^N} \sum_{(f,u) \in \text{data}} L(\hat{\A}(f; \theta), u),
\end{equation}
where $L$ is a loss function that measures the discrepancy between $\hat{A}(f;\theta)$ and $u$, and the sum is over all available training data pairs $(f,u)$. The challenges of operator learning often arise from selecting an appropriate neural operator architecture for $\hat{A}$, the computational complexities of solving the optimization problem, and the ability to generalize to new data.

A typical application of operator learning arises when learning the solution operator associated with a PDE, which maps a forcing function $f$ to a solution $u$. One can informally think of it as the (right) inverse of a differential operator. One of the simplest examples is the solution operator associated with Poisson's equation with zero Dirichlet conditions:
\begin{equation} \label{eq_poisson}
    -\nabla^2 u = f, \quad x\in \Omega\subset \R^d, \qquad u|_{\partial \Omega} = 0,
\end{equation}
where $\partial \Omega$ means the boundary of $\Omega$.  In this case, the solution operator, $\A$, can be expressed as an integral operator:
\[
    \A(f) = \int_{\Omega} G(\cdot,y)f(y)\d y = u,
\]
where $G$ is the Green's function associated with \cref{eq_poisson}~\citep[Chapt.~2]{evans2010partial}. A neural operator is then trained to approximate the action of $\A$ using training data pairs $(f_1,u_1),\ldots,(f_M,u_M)$.

In general, recovering the solution operator is challenging, as it is often nonlinear and high-dimensional, and the available data may be scarce or noisy. Nevertheless, unlike inverse problems, which aim to recover source terms from solutions, the forward problem is usually well-posed. As we shall see, learning solution operators lead to new insights or applications that can complement inverse problem techniques, as described in two surveys~\citep{stuart2010inverse} and \citep{arridge2019solving}.

\subsection{What is a neural operator?}

Neural operators~\citep{kovachki2023neural,lu2021learning} are analogues of neural networks with infinite-dimensional inputs. Neural operators were introduced to generalize standard deep learning techniques to learn mappings between function spaces instead of between discrete vector spaces $\R^{d_1}$ to $\R^{d_L}$, where $d_1$ is the input dimension of a neural network and $d_L$ is the output dimension. In its most traditional formulation, a fully connected neural network can be written as a succession of affine transformations and nonlinear activation functions as
\[\mathcal{N}(x) = \sigma(A_L(\cdots \sigma(A_1 x +b_1)\cdots) +b_L),\]
where $L\geq 1$ is the number of layers, $A_i$ are the weight matrices, $b_i$ are the bias vectors, and $\sigma:\R\to\R$ is the activation function, often chosen to be the ReLU function $\sigma(x)=\max\{x,0\}$. Neural operators generalize this architecture, where the input and output of the neural network are functions instead of vectors. Hence, the input of a neural operator is a function $f:\Omega\to \R^{d_1}$, where $\Omega\subset \R^d$ is the domain of the function, and the output is a function $u:\Omega\to \R^{d_L}$. The neural operator is then defined as a composition of integral operators and nonlinear functions, which results in the following recursive definition at layer $i$:
\begin{equation} \label{eq_layer_update}
    u_{i+1}(x) = \sigma\left(\int_{\Omega_i}K^{(i)}(x,y)u_i(y)\d y+b_i(x)\right),\quad x\in \Omega_{i+1},
\end{equation}
where $\Omega_i\subset\R^{d_i}$ is a compact domain, $b_i$ is a bias function, and $K^{(i)}$ is the kernel. The kernels and biases are then parameterized and trained similarly to standard neural networks. However, approximating the kernels or evaluating the integral operators could be computationally expensive. Hence, several neural operator architectures have been proposed to overcome these challenges, such as DeepONets~\citep{lu2021learning} and Fourier neural operators~\citep{li2020fourier}.

\subsection{Where is operator learning relevant?}

Operator learning has been successfully applied to many PDEs from different fields, including fluid dynamics with simulations of fluid flow turbulence in the Navier--Stokes equations at high Reynolds number~\citep{li2023long,peng2022attention}, continuum mechanics~\citep{you2022learning}, astrophysics~\citep{mao2023ppdonet}, quantum mechanics with the Schr\"odinger equation~\citep{li2020fourier}, and weather forecasting~\citep{kurth2023fourcastnet,lam2023learning}. The following four types of applications might directly benefit from operator learning.

\paragraph{Speeding up numerical PDE solvers.}
First, one can use operator learning to build reduced-order models of complex systems that are computationally challenging to simulate with traditional numerical PDE solvers. For example, this situation arises in fluid dynamics applications such as modeling turbulent flows, which require a very fine discretization or the simulation of high dimensional PDEs. Moreover, specific problems in engineering require the evaluation of the solution operator many times, such as in the design of aircraft or wind turbines. In these cases, a fast but less accurate solver provided by operator learning may be used for forecasting or optimization. This is one of the main motivations behind Fourier neural operators in~\citep{li2020fourier}. There are also applications of operator learning~\citep{zheng2023fast} to speed up the sampling process in diffusion models or score-based generative models~\citep{sohl2015deep,ho2020denoising,song2021scorebased}, which require solving complex differential equations. However, one must be careful when comparing performance against classical numerical PDE solvers, mainly due to the significant training time required by operator learning.

\paragraph{Parameter optimization.}
In our experience, the computational efficiency of operator learning is mainly seen in downstream applications such as parameter optimization. Once the solution operator has been approximated, it can be exploited in an inverse problem framework to recover unknown parameters of the PDE, which may be computationally challenging to perform with existing numerical PDE solvers. Additionally, neural operators do not rely on a fixed discretization as they are mesh-free and parameterized by a neural network that can be evaluated at any point. This property makes them suitable for solving PDEs on irregular domains or transferring the model to other spatial resolutions~\citep{kovachki2023neural}.

\paragraph{Benchmarking new techniques.}
Operator learning may also be used to benchmark and develop new deep learning models. As an example, one can design specific neural network architectures to preserve quantities of interest in PDEs, such as symmetries~\citep{olver1993applications}, conservation laws~\citep[Sec.~3.4]{evans2010partial}, and discretization independence. This could lead to efficient architectures that are more interpretable and generalize better to unseen data, and exploit geometric priors within datasets~\citep{bronstein2021geometric}. Moreover, the vast literature on PDEs and numerical solvers can be leveraged to create datasets and assess the performance of these models in various settings without requiring significant computational resources for training.

\paragraph{Discovering unknown physics.}
Last but not least, operator learning is helpful for the discovery of new physics~\citep{lu2021learning}. Indeed, the solution operator of a PDE is often unknown, or one may only have access to a few data points without any prior knowledge of the underlying PDE. In this case, operator learning can be used to discover the PDE or a mathematical model to perform predictions. This can lead to new insights into the system's behavior, such as finding conservation laws, symmetries, shock locations, and singularities~\citep{boulle2021data}. However, the complex nature of neural networks makes them challenging to interpret or explain, and there are many future directions for making SciML more interpretative.

\subsection{Organization of the paper} \label{sec_organization}

This paper is organized as follows. We begin in \cref{sec_numerical} by exploring the connections between numerical linear algebra and operator learning. Then, we review the main neural network architectures used to approximate operators in \cref{sec_architectures}. In \cref{sec_learning}, we focus on the data acquisition process, a crucial step in operator learning. We discuss the choice of the distribution of source terms used to probe the system, the numerical PDE solver, and the number of training data. Along the way, we analyze the optimization pipeline, including the possible choices of loss functions, optimization algorithms, and assessment of the results. Finally, in \cref{sec_physics}, we conclude with a discussion of the remaining challenges in the field that include the development of open-source software and datasets, the theoretical understanding of the optimization procedure, and the discovery of physical properties in operator learning, such as symmetries and conservation laws.

\section{From numerical linear algebra to operator learning} \label{sec_numerical}

There is a strong connection between operator learning and the recovery of structured matrices from matrix-vector products.  Suppose one aims to approximate the solution operator associated with a linear PDE using a single layer neural operator in the form of \cref{eq_layer_update} without the nonlinear activation function. Then, after discretizing the integral operator using a quadrature rule, it can be written as a matrix-vector product, where the integral kernel $K:\Omega\times \Omega\to\R$ is approximated by a matrix $A\in \R^{N\times N}$. Moreover, the structure of the matrix is inherited from the properties of the Green's function (see~\cref{tab_property_matrices}). The matrix's underlying structure---whether it is low-rank, circulant, banded, or hierarchical low-rank---plays a crucial role in determining the efficiency and approach of the recovery process. This section describes the matrix recovery problem as a helpful way to gain intuition about operator learning and the design of neural operator architectures (see~\cref{sec_architectures}). Another motivation for recovering structured solution operators is to ensure that the neural operators are fast to evaluate, which is essential in applications involving parameter optimization and benchmarking~\citep{kovachki2023neural}.

\begin{table}[htbp]
    \centering
    \begin{tabular}{p{5cm}cc}
        \hline
        Property of the PDE                      & Solution operator's kernel & Matrix structure \\
        \hline
        \hspace{1.5cm}\textbf{---}\hspace{1.5cm} & Globally smooth            & Low rank         \\
        Periodic BCs \& const.~coeffs            & Convolution kernel         & Circulant        \\
        Localized behavior                       & Off-diagonal decay         & Banded           \\
        Elliptic / Parabolic                     & Off-diagonal low rank      & Hierarchical     \\
        \hline
    \end{tabular}
    \caption{Solution operators associated with linear PDEs can often be represented as integral operators with a kernel called a Green's function. The properties of a linear PDE induce different structures on the Green's function, such as translation-invariant or off-diagonal low-rank. When these integral operators are discretized, one forms a matrix-vector product, and hence the matrix recovery problem can be viewed as a discrete analogue of operator learning. The dash in the first column and row means no PDE has a solution operator with a globally smooth kernel.}
    \label{tab_property_matrices}
\end{table}

Consider an unknown matrix $A\in\mathbb{R}^{N\times N}$ with a known structure such as rank-$k$ or banded. We assume that the matrix $A$ is a black box and cannot be seen, but that one can probe $A$ for information via matrix-vector products, i.e., the maps $x\mapsto Ax$ and $x\mapsto A^\top x$, with $A^\top$ representing the transpose of $A$. The matrix recovery problem is the task of approximating the matrix $A$ using as few queries to $x\mapsto Ax$ and $x\mapsto A^\top x$ as possible. Every matrix with $N$ columns can be deduced in a maximum of $N$ matrix-vector product queries, as $Ae_j$ for $1 \leq j \leq N$ returns the $j$th column, where $e_j$ denotes the $j$th standard basis vector. However, if the matrix $A$ has a specific structure such as low-rank, circulant, banded, or hierarchical low-rank, it is often possible to recover $A$ using far fewer queries. This section describes how to recover structured matrices efficiently using matrix-vector products. We prefer doing matrix recovery with Gaussian random vectors because the infinite-dimensional analogue of these vectors are random functions drawn from a Gaussian process, which is a widespread choice of training input data in operator learning.

\subsection{Low rank matrix recovery} \label{sec_low_rank}
Let $A\in \R^{N\times N}$ be a rank-$k$ matrix, then it can be expressed for some $C\in\mathbb{R}^{N\times k}$ and $R\in\mathbb{R}^{k\times N}$ as
\[
    A = CR.
\]
\citealt{halikias2022matrix} showed that at least $2k$ queries are required to capture the $k$-dimensional row and column spaces and deduce $A$. \citealt{halko2011finding,martinsson2020randomized} introduced the randomized singular value decomposition (SVD) as a method to recover a rank-$k$ matrix with probability one in $2k$ matrix-vector products with random Gaussian vectors. The randomized SVD can be expressed as a recovery algorithm in \cref{alg_rsvd}.

\renewcommand{\algorithmicrequire}{\textbf{Input:}}
\renewcommand{\algorithmicensure}{\textbf{Output:}}
\begin{algorithm}[htbp]
    \caption{Randomized singular value decomposition.}
    \label{alg_rsvd}
    \begin{algorithmic}[1]
        \State Draw a random matrix $X \in \mathbb{R}^{N\times k}$ with i.i.d. standard Gaussian entries.
        \State Perform $k$ queries with $A$: $Y = AX$.
        \State Compute the QR factorization of $Y = QR$.
        \State Perform $k$ queries with $A^\top$: $Z = A^\top Q$.
        \State Return $A = QZ^\top$.
    \end{algorithmic}
    \label{alg:rSVD}
\end{algorithm}

A randomized algorithm is crucial for low-rank matrix recovery to prevent input vectors from lying within the $N-k$ dimensional nullspace of $A$. Hence, recovering any low-rank matrix with a deterministic algorithm using fixed input vectors is impossible. Then, for the rank-$k$ matrix recovery problem,~\cref{alg:rSVD} recovers $A$ with probability one. A small oversampling parameter $p\geq 1$ is used for numerical stability, such as $p=5$. This means that $X \in \mathbb{R}^{N\times (k+p)}$, preventing the chance that a random Gaussian vector might be highly aligned with the nullspace of $A$.

A convenient feature of~\cref{alg:rSVD} is that it also works for matrices with numerical rank\footnote{For a fixed $0<\epsilon<1$, we say that a matrix $A$ has numerical rank $k$ if $\sigma_{k+1}(A)<\epsilon \sigma_1(A)$ and $\sigma_{k}(A)\geq \epsilon \sigma_1(A)$, where $\sigma_1(A)\geq \sigma_2(A)\geq\ldots\geq \sigma_N(A)\geq 0$ are the singular values of $A$.} $k$, provided that one uses a random matrix $X$ with $k+p$ columns. In particular, a simplified statement of~\citep[Thm.~10.7]{halko2011finding} shows that the randomized SVD recovers a near-best low-rank matrix in the sense that
\[
    \mathbb{P}\left[\|A - QZ^\top\|_{\textup{F}} \leq \left(1 + 15\sqrt{k+5}\right) \min_{{\rm rank}(A_k) \leq k} \|A - A_k\|_{\textup{F}}\right] \geq 0.999,
\]
where $\|\cdot\|_{\textup{F}}$ is the Frobenius norm of $A$. While other random embeddings can be used to probe $A$, Gaussian random vectors give the cleanest probability bounds~\citep{martinsson2020randomized}. Moreover, ensuring that the entries in each column of $X$ have some correlation and come from a multivariable Gaussian distribution allows for the infinite-dimensional extension of the randomized SVD and its application to recover Hilbert--Schmidt operators~\citep{boulle2022generalization,boulle2021learning}. This analysis allows one to adapt~\cref{alg:rSVD} to recover solution operators with low-rank kernels.

Low-rank matrix recovery is one of the most straightforward settings to motivate DeepONet (see~\cref{sec_deeponets}). One of the core features of DeepONet is to use the trunk net to represent the action of a solution operator on a set of basis functions generated by the so-called branch net. Whereas in low-rank matrix recovery, we often randomly draw the columns of $X$ as input vectors, DeepONet is trained with these functions. However, like DeepONet, low-rank matrix recovery is constructing an accurate approximant whose action is on these vectors. Many operators between function spaces can often be represented to high accuracy with DeepONet in the same way that the kernels of solution operators associated with linear PDEs often have algebraically fast decaying singular values.

\subsection{Circulant matrix recovery} \label{sec_circulant}

The FNO~\citep{li2020fourier} structure is closely related to circulant matrix recovery. Consider an $N \times N$ circulant matrix $C_c$, which is parameterized a vector $c\in \mathbb{R}^N$ as follows:
\[
    C_c =
    \begin{bmatrix}
        c_0     & c_{N-1} & \cdots  & c_2    & c_1     \\
        c_1     & c_0     & c_{N-1} & \ddots & c_2     \\
        \vdots  & c_1     & c_0     & \ddots & \vdots  \\
        c_{N-2} & \ddots  & \ddots  & \ddots & c_{N-1} \\
        c_{N-1} & c_{N-2} & \cdots  & c_1    & c_0
    \end{bmatrix}.
\]
To recover $C_c$ with a random Gaussian vector $g$, we recall that $C_c$ can be interpreted as a multiplication operator in the Fourier basis. By associativity of multiplication, we have
\[
    C_c g = C_g c.
\]
If we perform the matrix-vector product query $y=C_cg$, we can find the vector $c$ by solving the linear system $C_g c = y$.  Since $c$ completely defines $C_c$, we have recovered the circulant matrix. Moreover, the linear system $C_gc = y$ can be solved efficiently using the fast Fourier transform (FFT) in $\mathcal{O}(N\log N)$ operations. A convenient feature of circulant matrices is that given a new vector $x\in\mathbb{R}^N$, one can compute $C_cx$ in $\mathcal{O}(N\log N)$ operations using the FFT.

Circulant matrix recovery motivates Fourier neural operators. Hence, FNOs leverage the fast Fourier transform to efficiently parameterize the kernel of a solution operator, essentially capturing the operator in a spectral sense. Similarly, circulant matrices are diagonalized by the discrete Fourier transform matrix. The infinite-dimensional analogue of a circulant matrix is a solution operator with a periodic and translation invariant kernel, and this is the class of solution operators for which the FNO assumptions are fully justified. FNOs are extremely fast to evaluate because of their structure, making them popular for parameter optimization and favorable for benchmarking against reduced-order models.

\subsection{Banded matrix recovery} \label{sec_banded}

We now consider a banded matrix $A\in \R^{N\times N}$ with a fixed bandwidth $w$, i.e.,
\[
    A_{ij} = 0, \quad \text{if } |i-j|>w.
\]
The matrix $A$ can be recovered with $w+2$ matrix-vector products, but not fewer, using the $w+2$ columns of the following matrix as input vectors:
\[
    \begin{bmatrix}I_{w+2} & \cdots & I_{w+2} \end{bmatrix}^\top,
\]
where $I_{w+2}$ is the $(w+2)\times (w+2)$ identity matrix. Since every $w$th column has disjoint support, these input vectors recover the columns of $A$. Of course, this also means that a $A$ can be recovered with $w+2$ Gaussian random vectors.

\begin{figure}[htbp]
    \centering
    \scalebox{0.7}{
        \begin{minipage}[l]{0.49\textwidth}
            \[
                \left[
                    \begin{array}{cccccccccccc}
                        \color{myred}* & \color{mygreen} * & \color{myblue}* &                   &                   &                 &                  &                 &                   &                   &                &                  \\
                        \color{myred}* & \color{mygreen}*  & \color{myblue}* & \color{myyellow}* &                   &                 &                  &                 &                   &                   &                &                  \\
                        \color{myred}* & \color{mygreen}*  & \color{myblue}* & \color{myyellow}* & \color{myorange}* &                 &                  &                 &                   &                   &                &                  \\
                                       & \color{mygreen}*  & \color{myblue}* & \color{myyellow}* & \color{myorange}* & \color{myred} * &                  &                 &                   &                   &                &                  \\
                                       &                   & \color{myblue}* & \color{myyellow}* & \color{myorange}* & \color{myred}*  & \color{mygreen}*                                                                                               \\
                                       &                   &                 & \color{myyellow}* & \color{myorange}* & \color{myred}*  & \color{mygreen}* & \color{myblue}*                                                                             \\
                                       &                   &                 &                   & \color{myorange}* & \color{myred}*  & \color{mygreen}* & \color{myblue}* & \color{myyellow}*                                                         \\
                                       &                   &                 &                   &                   & \color{myred}*  & \color{mygreen}* & \color{myblue}* & \color{myyellow}* & \color{myorange}*                                     \\
                                       &                   &                 &                   &                   &                 & \color{mygreen}* & \color{myblue}* & \color{myyellow}* & \color{myorange}* & \color{myred}*                    \\
                                       &                   &                 &                   &                   &                 &                  & \color{myblue}* & \color{myyellow}* & \color{myorange}* & \color{myred}* & \color{mygreen}* \\
                                       &                   &                 &                   &                   &                 &                  &                 & \color{myyellow}* & \color{myorange}* & \color{myred}* & \color{mygreen}* \\
                                       &                   &                 &                   &                   &                 &                  &                 &                   & \color{myorange}* & \color{myred}* & \color{mygreen}*
                    \end{array}
                    \right]
            \]
        \end{minipage}
    }
    \put(-180,40){(a)}
    \put(5,40){(b)}
    \begin{minipage}[r]{0.49\textwidth}
        \begin{overpic}[width=\textwidth]{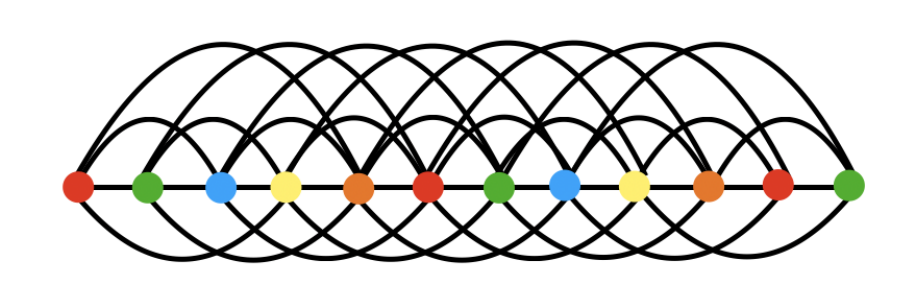}
        \end{overpic}
    \end{minipage}
    \caption{(a) A generic $12\times 12$ banded matrix with bandwidth $2$, with a maximum of $5$ diagonals, and the corresponding graph (b). Here, each vertex is a column of the banded matrix, and two vertices are connected if their corresponding columns do not have disjoint support. The coloring number of $5$ determines the minimum number of matrix-vector products needed to recover the structure. Generally, an $N\times N$ banded matrix with bandwidth $w$ can be recovered in $2w+1$ matrix-vector products.}
    \label{fig:bandedMatrixRecovery}
\end{figure}

There is a way to understand how many queries one needs as a graph-coloring problem. Consider the graph of size $N$, corresponding to an $N\times N$ banded matrix with bandwidth $w$, where two vertices are connected if their corresponding columns do not have disjoint support (see~\cref{fig:bandedMatrixRecovery}). Then, the minimum number of matrix-vector product queries needed to recover $A$ is the graph coloring number of this graph.\footnote{Recall that the coloring number of a graph is the minimum number of colors required to color the vertices so that no two vertices connected by an edge are identically colored.} One can see why this is the case because all the columns with the same color can be deduced simultaneously with a single matrix-vector product as they must have disjoint support.

Banded matrix recovery motivates Graph neural operators (GNOs), which we will describe later in \cref{sec_local}, as both techniques exploit localized structures within data. GNOs use the idea that relationships in nature are local and can be represented as graphs with no faraway connections. By only allowing local connections, GNOs can efficiently represent solution operators corresponding to local solution operators, mirroring the way banded matrices are concentrated on the diagonal. Likewise, with a strong locality, GNOs are relatively fast to evaluate, making them useful for parameter optimization and benchmarking. However, they may underperform if the bandwidth increases or the solution operator is not local.

\subsection{Hierarchical low rank matrix recovery} \label{sec_hierarchical}
\begin{figure}[ht]
    \centering
    \scalebox{0.68}{
        \begin{minipage}[l]{0.49\textwidth}
            \[
                \renewcommand{\arraystretch}{.2}
                \setlength{\arraycolsep}{.1pt}
                \begin{array}{|c|c|}
                    \hline
                    \begin{array}{c|c}
                        \begin{array}{c|c} \rule{0pt}{1.45\normalbaselineskip}  H_{11} & \ W_3Z_3^\top \\[10pt]\hline\rule{0pt}{1.45\normalbaselineskip} \ U_3V_3^\top & H_{22}\\[10pt]
                        \end{array} & W_{1}Z_{1}^\top \\[10pt]\hline\rule{0pt}{1.45\normalbaselineskip} U_{1}V_{1}^\top &
                        \begin{array}{c|c} \rule{0pt}{1.45\normalbaselineskip}H_{33} & \ W_4Z_4^\top \\[10pt]\hline\rule{0pt}{1.45\normalbaselineskip} \ U_4V_4^\top & H_{44}\\[10pt]
                        \end{array}
                    \end{array} & W_0Z_0^\top \\[10pt]\hline\rule{0pt}{1.45\normalbaselineskip} U_0V_0^\top &
                    \begin{array}{c|c}\rule{0pt}{1.45\normalbaselineskip}
                        \begin{array}{c|c} \rule{0pt}{1.5\normalbaselineskip}H_{55} & \ W_5Z_5^\top \\[10pt]\hline\rule{0pt}{1.45\normalbaselineskip} \ U_5V_5^\top & H_{66} \\[10pt]
                        \end{array} & W_{2}Z_{2}^\top \\[10pt]\hline\rule{0pt}{1.45\normalbaselineskip} U_{2}V_{2}^\top &
                        \begin{array}{c|c} \rule{0pt}{1.45\normalbaselineskip}H_{77} & \ W_6Z_6^\top \\[10pt]\hline\rule{0pt}{1.45\normalbaselineskip} \ U_6V_6^\top & H_{88}\\[10pt]
                        \end{array} \\
                    \end{array} \\
                    \hline
                \end{array}
            \]
        \end{minipage}}
    \put(-165,80){(a)}
    \hspace{2cm}
    \begin{minipage}[r]{0.44\textwidth}
        \vspace{0.2cm}
        \begin{overpic}[width=\textwidth,trim=0 0 0 0,clip]{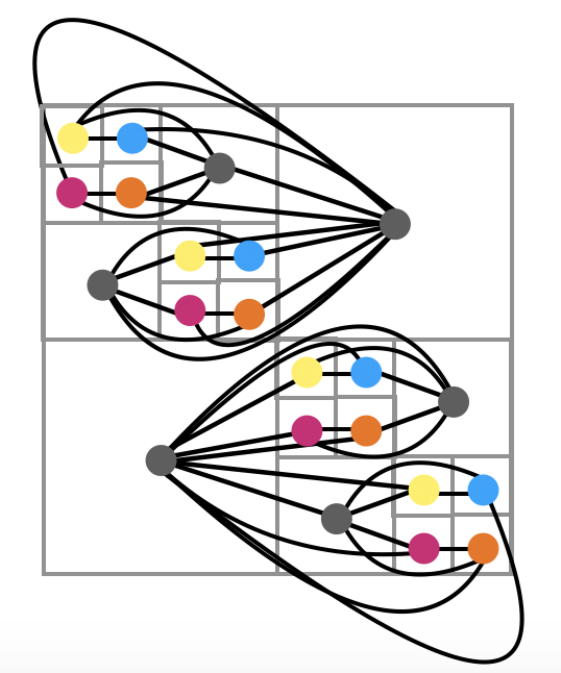}
        \end{overpic}
    \end{minipage}
    \put(-205,80){(b)}
    \caption{(a) A HODLR matrix $H_{N,k}$ after three levels of partitioning. Since $H_{N,k}$ is a rank-$k$ HODLR matrix, $U_i$, $V_i$, $W_i$, and $Z_i$ have at most $k$ columns. The matrices $A_{ii}$ are themselves rank-$k$ HODLR matrices of size $N/8\times N/8$ and can be further partitioned. (b) Graph corresponding to a hierarchical low-rank matrix with three levels. Here, each vertex is a low-rank block of the matrix, where two vertices are connected if their low-rank blocks occupy the same row. At each level, the number of required matrix-vector input probes to recover that level is proportional to the coloring number of the graph when restricted to submatrices of the same size. In this case, the submatrices that are identically colored can be recovered simultaneously.}
    \label{fig:StructureDiagram}
\end{figure}

An $N \times N$ rank-$k$ hierarchical off-diagonal low rank (HODLR) matrix, denoted as $H_{N,k}$, is a structure that frequently appears in the context of discretized solution operators associated with elliptic and parabolic PDEs~\citep{hackbusch2004hierarchical}. To understand its recursive structure, we assume $N$ to be a power of 2 and illustrate the structure in~\cref{fig:StructureDiagram}(a).

Since a rank-$k$ matrix requires $2k$ matrix-vector product queries to be recovered (see \cref{sec_low_rank}), a naive approach to deducing $H_{N,k}$ is to use $2k$ independent queries on each submatrix. However, one can show that some of the submatrices of $H_{N,k}$ can be recovered concurrently using the same queries. We use a graph coloring approach~\citep{levitt2022randomized} to determine which submatrices can be recovered concurrently. This time, we consider the graph where each vertex is a low-rank submatrix of $H_{N,k}$ and connect two vertices if their corresponding low-rank submatrices occupy the same column as in~\cref{fig:StructureDiagram}(b). The low-rank submatrices that are identically colored at each level can be recovered concurrently in only $2k$ queries. Hence, it can be shown that an $N\times N$ hierarchical rank-$k$ matrix can be recovered in fewer than $10k\lceil \log_2(N)\rceil$ matrix-vector products~\citep{halikias2022matrix}.  The precise coloring of the graph in~\cref{fig:StructureDiagram}(b) can also be used to derive a particular algorithm for hierarchical matrix recovery known as peeling~\citep{levitt2022linear,levitt2022randomized,lin2011fast,martinsson2011fast}. These peeling algorithms have been recently generalized to the infinite-dimensional setting by~\citep{boulle2023elliptic}.

HODLR recovery can be seen as the simplest version of a multipole graph neural operator (MGNO) (see~\cref{sec_multipole}) as both emphasize the importance of capturing operators at multiple scales. MGNOs are based on a hierarchical graph with interactions at different scales or levels (see~\cref{sec_local}). By incorporating local (near-field) and global (far-field) interactions, MGNOs can effectively learn complex patterns. MGNOs are often great at representing solution operators due to their multiscale nature. The price to pay is that the final neural operator can be computationally expensive to evaluate, and it is a complicated structure to implement.

\section{Neural operator architectures} \label{sec_architectures}

In this section, we review the main neural operator architectures used in the literature, namely DeepONets~\citep{lu2021learning}, Fourier neural operators~\citep{li2020fourier}, and Deep Green networks~\citep{gin2020deepgreen,boulle2021data}. We also refer to the recent survey by~\citep{goswami2023physics} for a review of the different neural operator architectures and their applications. Each of these architectures employs different discretization and approximation techniques to make the neural operator more efficient and scalable by enforcing certain structures on the kernel such as low-rank, periodicity, translation invariance, or hierarchical low-rank structure (see~\cref{tab_property_kernels}).

\begin{table}[htbp]
    \centering
    \begin{tabular}{ p{4cm} cc  }
        \hline
        Neural operators      & Property of the operator & Kernel parameterization   \\
        \hline
        DeepONet              & Low-rank                 & Branch and trunk networks \\
        FNO                   & Translation-invariant    & Fourier coefficients      \\
        GreenLearning         & Linear                   & Rational neural network   \\
        DeepGreen             & Semi-linear              & Kernel matrix             \\
        Graph neural operator & Diagonally dominant      & Message passing network   \\
        Multipole GNO         & Off-diagonal low rank    & Neural network            \\
        \hline
    \end{tabular}
    \caption{Summary table of neural operator architectures, describing the property assumption on the operator along with the discretization of the integral kernels.}
    \label{tab_property_kernels}
\end{table}

Most neural operator architectures also come with theoretical guarantees on their approximation power.  These theoretical results essentially consist of universal approximation properties for neural operators~\citep{chen1995universal,kovachki2023neural,lu2021learning}, in a similar manner as neural networks~\citep{devore1998nonlinear}, and quantitative error bounds based on approximation theory to estimate the size, i.e., the number of trainable parameters, of a neural operator needed to approximate a given operator between Banach spaces to within a prescribed accuracy~\citep{lanthaler2021error,yarotsky2017error}.

\subsection{Deep operator networks} \label{sec_deeponets}

Deep Operator Networks (DeepONets) are a promising model for learning nonlinear operators and capturing the inherent relationships between input and output functions~\citep{lu2021learning}. They extend the capabilities of traditional deep learning techniques by leveraging the expressive power of neural networks to approximate operators in differential equations, integral equations, or more broadly, any functional maps from one function space to another. A key theoretical motivation for DeepONet is the universal operator approximation theorem~\citep{chen1995universal,lu2021learning}. This result can be seen as an infinite dimensional analogue of the universal approximation operator for neural networks~\citep{cybenko1989approximation,hornik1991approximation}, which guarantee that a sufficiently wide neural network can approximate any continuous function to any accuracy. Since the introduction of DeepONets by~\citep{lu2021learning}, several research works focused on deriving error bounds for the approximation of nonlinear operators by DeepONets in various settings, such as learning the solution operator associated with Burger's equation or the advection-diffusion equation~\citep{deng2022approximation}, and the approximation of nonlinear parabolic PDEs~\citep{de2022generic,lanthaler2021error}.

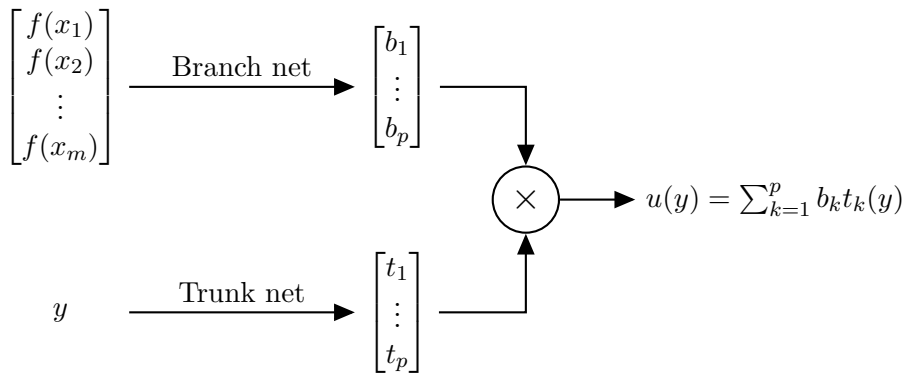
\begin{figure}[htbp]
    \centering
    \begin{tikzpicture}[auto, thick, >=triangle 45, on grid=false]
    \draw
    node at (0,0) [name=input1] {$\begin{bmatrix}f(x_1)\\f(x_2)\\\vdots\\f(x_m)\end{bmatrix}$}
    node [below of = input1, node distance = 3cm] (input2) {$y$}
    node [right = 3cm of input1] (suma1) {$\begin{bmatrix}b_1\\\vdots\\b_p\end{bmatrix}$}
    node [below of = suma1, node distance = 3cm] (suma2) {$\begin{bmatrix}t_1\\\vdots\\t_p\end{bmatrix}$}
    node [right = 1cm of suma1] (inte1) {}
    node [sum, below of = inte1, node distance = 1.5cm] (suma4) {\muma}
    node [right = 1cm of suma4] (ret2) {$u(y)=\sum_{k=1}^p b_k t_k(y)$};

    \draw[->](input1) -- node {Branch net}(suma1);
    \draw[->, transform canvas={yshift=-3cm}] (input1) -- (suma1) node[midway, above] {Trunk net};

    \draw[->](suma4) -- (ret2);

    \draw[->] (suma1) -| (suma4);
    \draw[->] (suma2) -| (suma4);

\end{tikzpicture}
    \caption{Schematic diagram of a deep operator network (DeepONet). A DeepONet parametrizes a neural operator using a branch network and a truncation (trunk) network. The branch network encodes the input function $f$ as a vector of $p$ features, which is then multiplied by the trunk network to yield a rank-$p$ representation of the solution $u$.}
    \label{fig_don}
\end{figure}

A DeepONet is a two-part deep learning network consisting of a branch network and a trunk network. The branch net encodes the operator's input functions $f$ into compact, fixed-size latent vectors $b_1(f(x_1),\ldots,f(x_m)),\ldots,b_p(f(x_1),\ldots,f(x_m))$, where $\{x_i\}_{i=1}^m$ are the sensor points at which the input functions are evaluated. The trunk net decodes these latent vectors to produce the final output function at the location $y\in \Omega$ as
\[\mathcal{N}(f)(y) = \sum_{k=1}^p b_k(f(x_1),\ldots,f(x_m)) t_k(y).\]
A schematic of a deep operator network is given in \cref{fig_don}.
The defining feature of DeepONets is their ability to handle functional input and output, thus enabling them to learn a wide array of mathematical operators effectively. It's worth mentioning that the branch network and the trunk network can have distinct neural network architectures tailored for different purposes, such as performing a feature expansion on the input of the trunk network as
$y\to \begin{pmatrix}
        y & \cos(\pi y) & \sin(\pi y) & \ldots
    \end{pmatrix}$
to take into account any potential oscillatory patterns in the data~\citep{di2023neural}. Moreover, while the interplay of the branch and trunk networks is crucial, the output of a DeepONet does not necessarily depend on the specific input points but rather on the global property of the entire input function, which makes it suitable for learning operator maps.

One reason behind the performance of DeepONet might be its connection with the low-rank approximation of operators and the SVD (see \cref{sec_low_rank}). Hence, one can view the trunk network as learning a basis of functions $\{t_k\}_{k=1}^p$ that are used to approximate the operator, while the branch network expresses the output function in this basis by learning the coefficients $\{b_k\}_{k=1}^p$. Moreover, the branch network can be seen as a feature extractor, which encodes the input function into a compact representation, thus reducing the problem's dimensionality to $p$, where $p$ is the number of branch networks. Additionally, several architectures, namely the POD-DeepONet~\citep{lu2022comprehensive} and SVD-DeepONet~\citep{venturi2023svd}, have been proposed to strengthen the connections between DeepONet and the SVD of the operator and increase its interpretability.

A desirable property for a neural operator architecture is to be discretization invariant in the sense that the model can act on any discretization of the source term and be evaluated at any point of the domain~\citep{kovachki2023neural}. This property is crucial for the generalization of the model to unseen data and the transferability of the model to other spatial resolutions. While DeepONets can be evaluated at any location of the output domain, DeepONets are not discretization invariant in their original formulation by~\citealt{lu2021learning} as the branch network is evaluated at specific points of the input domain (see~\cref{fig_don}). However, this can be resolved using a low-rank neural operator~\citep{kovachki2023neural}, sampling the input functions at local spatial averages~\citep{lanthaler2021error}, or employing a principal component analysis (PCA) alternative of the branch network~\citep{de2022cost}.

The training of DeepONets is performed using a supervised learning process. It involves minimizing the mean-squared error between the predicted output $\mathcal{N}(f)(y)$ and the actual output of $u$ the operator on the training functions at random locations $\{y_j\}_{j=1}^{n}$, i.e.,
\begin{equation} \label{eq_loss_don}
    \min_{\theta\in \R^N}\frac{1}{|\text{data}|}\sum_{(f,u)\in\text{data}}\frac{1}{n}\sum_{j=1}^n |\mathcal{N}(f)(y_j)-u(y_j)|^2.
\end{equation}
The term inside the first sum approximates the integral of the mean-squared error, $|\mathcal{N}(f)-u|^2$, over the domain $\Omega$ using Monte-Carlo integration. The optimization is typically done via backpropagation and gradient descent algorithms, which are the same as in traditional neural networks. Importantly, DeepONets allow for different choices of loss functions, depending on the problem. For example, mean squared error is commonly used for regression tasks, but other loss functions might be defined to act as a regularizer and incorporate prior physical knowledge of the problem~\citep{goswami2022physics,wang2021learning}. The selection of an appropriate loss function is a crucial step in defining the learning process of these networks and has a substantial impact on their performance (see~\cref{sec_loss}).

DeepONet has been successfully applied and adapted to a wide range of problems, including predicting cracks in fracture mechanics using a variational formulation of the governing equations~\citep{goswami2022physics}, simulating the New York-New England power grid behavior with a probabilistic and Bayesian framework to quantify the uncertainty of the trajectories~\citep{moya2023deeponet}, as well as predicting linear instabilities in high-speed compressible flows with boundary layers~\citep{di2023neural}.

\subsection{Fourier neural operators} \label{sec_FNO}

Fourier neural operators (FNOs)~\citep{li2020fourier,kovachki2023neural} are a class of neural operators motivated by Fourier spectral methods. FNOs have found their niche in dealing with high-dimensional PDEs, which are notoriously difficult to solve using traditional numerical methods due to the curse of dimensionality. They've demonstrated significant success in learning and predicting solutions to various PDEs, particularly those with periodic boundary conditions or those that can be transformed into the spectral domain via Fourier transform. This capability renders FNOs an invaluable tool in areas where PDEs play a central role, such as fluid dynamics, quantum mechanics, and electromagnetism.

The main idea behind FNOs is to choose the kernels $K^{(i)}$ in \cref{eq_layer_update} as translation-invariant kernels satisfying $K^{(i)}(x,y)=k^{(i)}(x-y)$ (provided the input and output domains are torus) such that the integration of the kernel can be performed efficiently as a convolution using the Fast Fourier Transform (FFT)~\citep{cooley1965algorithm}, i.e., multiplication in the feature space of Fourier coefficients. Hence, the integral operation in \cref{eq_layer_update} can be performed as
\[\int_{\Omega_i} k^{(i)}(x-y)u_i \d y=\mathcal{F}^{-1}(\mathcal{F}(k^{(i)})\mathcal{F}(u_i))(x)=\mathcal{F}^{-1}(\mathcal{R}\cdot\mathcal{F}(u_i))(x),\quad x\in \Omega_i,\]
where $\mathcal{F}$ denotes the Fourier transform and $\mathcal{F}^{-1}$ its inverse. The kernel $K^{(i)}$ is parametrized by a periodic function $k^{(i)}$, which is discretized by a (trainable) weight vector of Fourier coefficients $\mathcal{R}$, and truncated to a finite number of Fourier modes. Then, if the input domain is discretized uniformly with $m$ sensor points, and the vector $\mathcal{R}$ contains at most $k_{\max}\leq m$ modes, the convolution can be performed in quasi-linear complexity in $\mathcal{O}(m\log m)$ operations via the FFT. This is a significant improvement over the $\mathcal{O}(m^2)$ operations required to evaluate the integral in \cref{eq_layer_update} using a quadrature rule. In practice, one can restrict the number of Fourier modes to $k_{\max}\ll m$ without significantly affecting the accuracy of the approximation whenever the input and output functions are smooth so that their representation in the Fourier basis enjoy rapid decay of the coefficients, thus further reducing the computational and training complexity of the neural operator.

\begin{figure}[htbp]
    \centering
    \begin{tikzpicture}[auto, thick, >=triangle 45, on grid=false]
    \draw
    node at (0,0) [name=input1] {$f(x)$}
    node [block, right = 1cm of input1] (suma1) {$v_0$}
    node [right = 2.5cm of suma1] (inte1) {\ldots}
    node [block, right = 2.5cm of inte1] (ret1) {$v_T$}
    node [right = 1cm of ret1] (ret2) {$u(x)$};
    \node (a) at ($(suma1)!.5!(inte1)$) {};
    \draw[color=gray,thick] (a.center) -- (a |- 0,-1) node {};

    \draw[->](input1) -- node {$P$}(suma1);
    \draw[->](suma1) -- node {Fourier layer} (inte1);
    \draw[->](inte1) -- node {Fourier layer} (ret1);
    \draw[->](ret1) -- node {$Q$} (ret2);

    \draw
    node [block, below of = inte1, node distance = 2cm] (inte2) {$\mathcal{F}^{-1}(\mathcal{R}\cdot\mathcal{F}(v))+b$}
    node [block, below of = suma1, node distance = 2cm] [name=suma3] {$v$}
    node [sum, below of = ret1, node distance = 2cm] (suma4) {\suma}
    node [block, right = 1cm of suma4] (inte3) {$\sigma$}
    node [block, below of =inte2, name=ret2, node distance = 2cm] {$W$}
    ;
    \draw[->] (suma3) -- node {} (inte2);
    \draw[->] (inte2) -- node {} (suma4);
    \draw[->] (suma4) -- node {} (inte3);
    \draw[->] (suma3) |- (ret2);
    \draw[->] (ret2) -| (suma4);

    \draw [color=gray,thick](1,-5) rectangle (12.2,-1);
    \node at (2.5,-4.7) [] {Fourier layer};
\end{tikzpicture}
    \caption{Schematic diagram of a Fourier neural operator (FNO). The networks P and Q, respectively, lift the input function $f$ to a higher dimensional space and project the output of the last Fourier layer to the output dimension. An FNO mainly consists of a succession of Fourier layers, which perform the integral operations in neural operators as a convolution in the Fourier domain and component-wise composition with an activation function $\sigma$.}
    \label{fig_diag_FNO}
\end{figure}
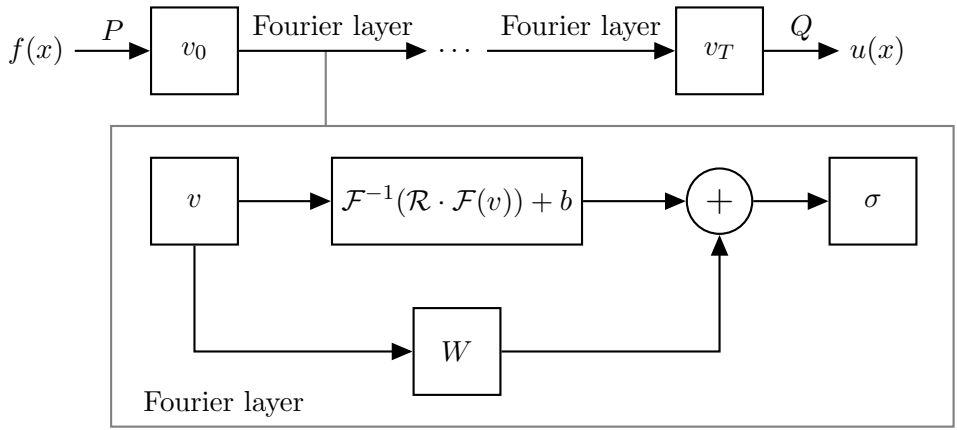

We display a diagram of the architecture of an FNO in \cref{fig_diag_FNO}. The input function $f$ is first lifted to a higher dimensional space by a neural network $P$. Then, a succession of Fourier layers is applied to the lifted function, which is parametrized by a vector of Fourier coefficients $\mathcal{R}_i$, a bias vector $b_i$, and a weight matrix $W$. Then, the output of the FNO at the $i$th layer is given by
\[v_{i} = \sigma(W_i v_{i-1}+\mathcal{F}^{-1}(\mathcal{R}_{i}\cdot\mathcal{F}(u_{i-1}))+b_i),\]
where $\sigma:\R\to\R$ is the activation function whose action is defined component-wise, often chosen to be the ReLU function. The weight matrix $W_i$ and bias vector $b_i$ perform a linear transformation of the input function $v_i$. After the last Fourier layer, the output of the FNO is obtained by applying a final neural network $Q$ on the output of the last Fourier layer to project it to the output dimension.

The training of FNOs, like DeepONets, is carried out via a supervised learning process. It typically involves minimizing a loss function that measures the discrepancy between the predicted and the true output of the operator on the input functions with respect to the trainable parameters of the neural network as in \cref{eq_loss_don}. Here, one needs to perform backpropagation through the Fourier layers, which is enabled by the implementation of fast GPU differentiable FFTs~\citep{mathieu2014fast} in deep learning frameworks such as PyTorch~\citep{paszke2019pytorch} and TensorFlow~\citep{abadi2016tensorflow}.

While FNOs have been proposed initially to alleviate the computational expense of performing integral operations in neural operators by leveraging the FFT, they have a distinctive advantage in learning operators where computations in the spectral domain are more efficient or desirable. This arises naturally when the target operator, along with the input and output functions, are smooth so that their representation as Fourier coefficients decay exponentially fast, yielding an efficient truncation. Hence, by selecting the architecture of the FNO appropriately, such as the number of Fourier modes $k_{\max}$, or the initialization of the Fourier coefficients, one can obtain a neural operator that preserves specific smoothness properties. However, when the input or output training data is not smooth, FNO might suffer from Runge's phenomenon near discontinuities~\citep{de2022cost}.

One main limitation of the FNO architectures is that the FFT should be performed on a uniform grid and rectangular domains, which is not always the case in practice. This can be overcome by applying embedding techniques to transform the input functions to a uniform grid and extend them to simple geometry, using a Fourier analytic continuation technique~\citep{bruno2007accurate}. Recently, several works have been proposed to extend the FNO architecture to more general domains, such as using a zero padding, linear interpolation~\citep{lu2022comprehensive}, or encoding the geometry to a regular latent space with a neural network~\citep{li2022fourier,li2023geometry}. However, this might lead to a loss of accuracy and additional computational cost. Moreover, the FFT is only efficient for approximating translation invariant kernels, which do not occur when learning solution operators of PDEs with non-constant coefficients.

Other related architectures aim to approximate neural operators directly in the feature space, such as spectral neural operators (SNO)~\citep{fanaskov2022spectral}, which are based on spectral methods and employ a simple feedforward neural network to map the input function, represented as a vector of Fourier or Chebyshev coefficients to an output vector of coefficients. Finally, \citep{raonic2023convolutional} introduces convolutional neural operators (CNOs) to alleviate the aliasing phenomenon of convolutional neural networks (CNNs) by learning the mapping between bandlimited functions. Contrary to FNOs, CNOs parameterize the integral kernel on a $k\times k$ grid and perform the convolution in the physical space as
\[\int_{\Omega}k(x-y)f(y)\d y = \sum_{i,j=1}^k k_{ij}f(x-z_ij),\quad x\in \Omega,\]
where $z_{ij}$ are the grid points.

Similarly to DeepONets, Fourier neural operators are universal approximators, in the sense that they are dense in the space of continuous operators~\citep{bhattacharya2021model,kovachki2023neural,kovachki2021universal}. However, even while being universal approximators, FNOs could, in theory, require a huge number of parameters to approximate a given operator to a prescribed accuracy $\epsilon>0$. As an example, \citep{kovachki2021universal} showed that the size of the FNO must grow exponentially fast as $\epsilon$ decreases to approximate any operator between rough functions whose Fourier coefficients decay only at a logarithmic rate. Fortunately, these pessimistic lower bounds are not observed in practice when learning solution operators associated with PDEs. Indeed, in this context, one can exploit PDE regularity theory and Sobolev embeddings to derive quantitative bounds on the size of FNOs for approximating solution operators that only grow sublinearly with the error. Here, we refer to the analysis of Darcy flow and the two-dimensional Navier--Stokes equations by~\citep{kovachki2021universal}.

\subsection{Deep Green networks} \label{sec_deep_green}

Deep Green networks (DGN) employ a different approach compared to DeepONets and FNOs to approximate solution operators of PDEs~\citep{gin2020deepgreen,boulle2021data}. Instead of enforcing certain properties on the integral kernel in \cref{eq_layer_update}, such as being low-rank (DeepONet) or translation-invariant (FNO), DGN learns the kernel directly in the physical space. Hence, assume that the underlying differential operator is a linear boundary value problem of the form
\begin{equation} \label{eq_integral_op}
    \Lop u = f \text{ in }\Omega,\quad\text{and}\quad u = 0 \text{ on }\partial\Omega.
\end{equation}
Under suitable regularity assumptions on operator $\Lop$ (e.g., uniform ellipticity or parabolicity), the solution operator $\A$  can be expressed as an integral operator with a Green kernel $G:\Omega\times\Omega\to \R\cup\{\infty\}$ as~\citep{evans2010partial,boulle2022learning,boulle2021learning}
\[
    \mathcal{A}(f)(x)=u(x)=\int_{\Omega} G(x,y)f(y)\d y,\quad x\in \Omega.
\]
\citep{boulle2021data} introduced GreenLearning networks (GL) to learn the Green kernel $G$ directly from data. The main idea behind GL is to parameterize the kernel $G$ as a neural network $\mathcal{N}$ and minimize the following relative mean-squared loss function to recover an approximant to $G$:
\begin{equation} \label{eq_loss_gl}
    \min_{\theta\in \R^N}\frac{1}{|\text{data}|}\sum_{(f,u)\in\text{data}}\frac{1}{\|u\|_{L^2(\Omega)}^2}\int_{\Omega}\left(u(x)-\int_\Omega \mathcal{N}(x,y)f(y)\d y\right)^2\d x.
\end{equation}
Once trained, the network $\mathcal{N}$ can be evaluated at any point in the domain, similarly to FNO and DON. A key advantage of this method is that it provides a more interpretable model, as the kernel can be visualized and analyzed to recover properties of the underlying differential operators~\citep{boulle2021data}. However, this comes at the cost of higher computational complexity, as the integral operation in \cref{eq_loss_gl} must be computed accurately using a quadrature rule and typically requires $\mathcal{O}(m^2)$ operations, as opposed to the $\mathcal{O}(m\log m)$ operations required by FNOs, where $m$ is the spatial discretization of the domain $\Omega$.

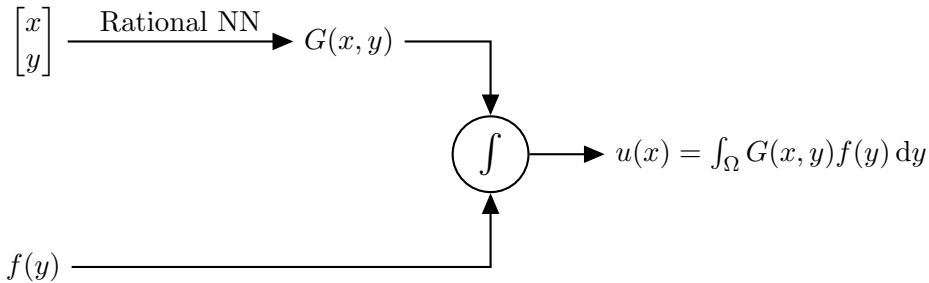
\begin{figure}[htbp]
    \centering
    \begin{tikzpicture}[auto, thick, >=triangle 45, on grid=false]
    \draw
    node at (0,0) [name=input1] {$\begin{bmatrix}x\\y\end{bmatrix}$}
    node [below of = input1, node distance = 3cm] (input2) {$f(y)$}
    node [right = 3cm of input1] (suma1) {$G(x,y)$}
    node [right = 1cm of suma1] (inte1) {}
    node [sum, below of = inte1, node distance = 1.5cm] (suma4) {\inta}
    node [right = 1cm of suma4] (ret2) {$u(x)=\int_\Omega G(x,y)f(y)\d y$};

    \draw[->](input1) -- node {Rational NN}(suma1);

    \draw[->](suma4) -- (ret2);

    \draw[->] (suma1) -| (suma4);
    \draw[->] (input2) -| (suma4);

\end{tikzpicture}
    \caption{Schematic diagram of a GreenLearning network (GL), which approximates the integral kernel (Green's function) associated with linear PDEs using rational neural networks. }
    \label{fig_gl}
\end{figure}

As Green's functions may be unbounded or singular, \citep{boulle2021data} propose to use a rational neural network~\citep{boulle2020rational} to approximate the Green kernel. A rational neural network is a neural network whose activation functions are rational functions, defined as the ratio of two polynomials whose coefficients are learned during the training phase of the network. This choice of architecture is motivated by the fact that rational networks have higher approximation power than the standard ReLU networks~\citep{boulle2020rational}, in the sense that they require exponentially fewer layers to approximate continuous functions within a given accuracy and may take arbitrary large values, which is a desirable property for approximating Green kernels. A schematic diagram of a GL architecture is available in \cref{fig_gl}.

When the underlying differential operator is nonlinear, the solution operator $\A$ cannot be written as an integral operator with a Green's function. In this case, \citep{gin2020deepgreen} propose to learn the solution operator $\A$ using a dual auto-encoder architecture Deep Green network (DGN), which is a neural operator architecture that learns an invertible coordinate transform map that linearizes the nonlinear boundary value problem. The resulting linear operator is approximated by a matrix, which represents a discretized Green's function but could also be represented by a neural network if combined with the GL technique. This approach has been successfully applied to learn the solution operator of the nonlinear cubic Helmholtz equation non-Sturm--Liouville equation and discover an underlying Green's function~\citep{gin2020deepgreen}.

Other deep learning-based approaches~\citep{lin2023bi,peng2023deep,sun2023binn} have since been introduced to recover Green's functions using deep learning, but they rely on a PINN technique in the sense that they require the knowledge of the underlying PDE operator. Finally, \citep{stepaniants2023learning} proposes to learn the Green kernel associated with linear partial differential operators using a reproducible kernel Hilbert space (RKHS) framework, which leads to a convex loss function.

\subsection{Graph neural operators} \label{sec_local}

As described in \cref{sec_deep_green}, solution operators associated with linear, elliptic, or parabolic PDEs of the form $\Lop u = f$ can be written as an integral operator with a Green kernel $G$~\citep[Sec.~2.2.4]{evans2010partial}. For simplicity, we consider Green kernels associated with uniformly elliptic operators in divergence form defined on a bounded domain $\Omega$ in spatial dimension $d\geq 3$:
\begin{equation} \label{eq_div_form}
    \Lop u = -\div(A(x)\nabla u) = f,\quad \text{on }\Omega\subset \R^d,
\end{equation}
where $A(x)$ is a bounded coefficient matrix satisfying the uniform ellipticity condition $A(x)\xi\cdot\xi\geq \lambda |\xi|^2$ for all $x\in \Omega$ and $\xi\in \R^d$, for some $\lambda>0$. In this section, we present a neural operator architecture that takes advantage of the local structure of the Green kernel associated with \cref{eq_div_form}, inferred by PDE regularity theory.

This architecture is called graph neural operator (GNO)~\citep{li2020neural} and is inspired by graph neural network (GNN) models~\citep{scarselli2008graph,wu2020comprehensive,zhou202057}. It focuses on capturing the Green kernel's short-range interactions to reduce the integral operation's computational complexity in \cref{eq_layer_update}. The main idea behind GNO is to perform the integral operation in \cref{eq_integral_op} locally on a small ball of radius $r$, $B(x,r)$, around $x$ for each $x\in\Omega$ as follows:
\begin{equation} \label{eq_local_integral}
    \mathcal{A}(f)(x)=u(x)\approx \int_{B(x,r)} G(x,y)f(y)\d y,\quad x\in \Omega.
\end{equation}
Here, \citep{li2020neural} propose to discretize the domain $\Omega$ using a graph, whose nodes represent discretized spatial locations, and use a message passing network architecture~\citep{gilmer2017neural} to perform an average aggregation of the nodes as in \cref{eq_local_integral}. The approach introduced by~\citep{li2020neural} aims to approximate the restriction $G_r$ to the Green's function $G$ on a band of radius $r$ along the diagonal of the domain $\Omega\times \Omega$, defined as
\[
    G_r(x,y) = \begin{cases}
        G(x,y), & \text{if }|x-y|\leq r, \\
        0,      & \text{otherwise},
    \end{cases}
\]
where $|\cdot|$ is the Euclidean distance in $\R^d$.

\begin{figure}[htbp]
    \centering
    \begin{overpic}[width=\textwidth]{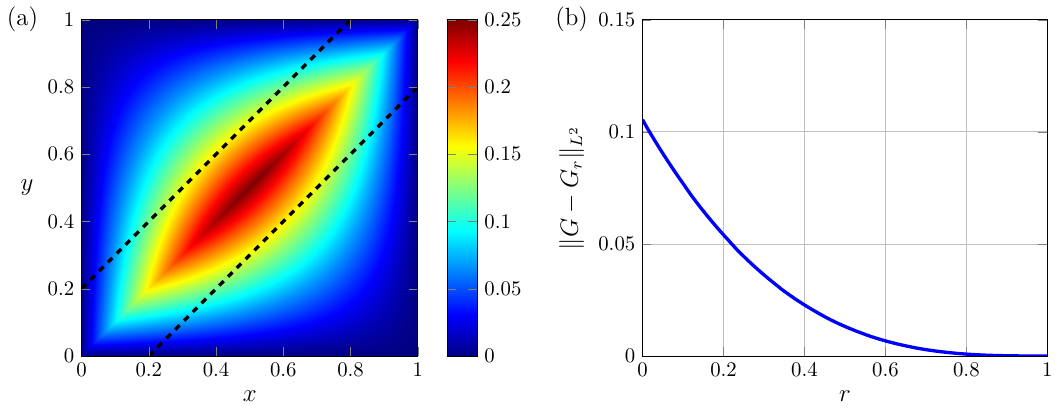}
    \end{overpic}
    \caption{(a) Green's function associated with the one-dimensional Poisson equation on $\Omega=[0,1]$ with homogeneous Dirichlet boundary conditions. The dashed lines highlight a band of radius $r=\sqrt{2}/10$ around the diagonal. (b) $L^2$-norm of the error between the Green's function $G$ and its truncation $G_r$ along a bandwidth of radius $r$ along the diagonal of the domain.}
    \label{fig_green_diagonal}
\end{figure}

This neural architecture is justified by the following pointwise bound satisfied by the Green's function and proven by~\citet[Thm.~1.1]{gruter1982green}:
\begin{equation} \label{eq_green_bound}
    |G(x,y)|\leq C(d,A)|x-y|^{2-d},\quad x,y\in \Omega,
\end{equation}
where $\Omega$ is a compact domain in $\R^d$ for $d\geq 3$, and $C$ is a constant depending only on $d$ and the coefficient matrix $A(x)$. Similar bounds have been derived in spatial dimension $d=2$ by~\citep{cho2012global,dong2009green} and for Green's functions associated with time-dependent, parabolic, PDEs~\citep{hofmann2004gaussian,cho2008green}. Then, integrating \cref{eq_green_bound} over the domain $\Gamma_r\coloneqq\{(x,y)\in \Omega\times\Omega\colon|x-y|>r\}$ yields a bound on the approximation error between $G$ and $G_r$ that decays algebraically fast as $r$ increases:
\begin{align*}
    \|G-G_r\|_{L^2(\Omega\times\Omega)} & =  \left(\int_{\Gamma_r}|G(x,y)|^2\d x\d y\right)^{1/2}
    \leq C(d,A)\left(\int_{\Gamma_r}r^{4-2d}\d x\d y\right)^{1/2}                                 \\
                                        & \leq |\Omega|C(d,A)r^{2-d}.
\end{align*}
This implies that the Green's function can be well approximated by a bandlimited kernel $G_r$ and that the approximation error bound improves in high dimensions. To illustrate this, we plot in \cref{fig_green_diagonal} the Green's function associated with the one-dimensional Poisson equation on $\Omega=[0,1]$ with homogeneous Dirichlet boundary conditions, along with the error between the Green's function $G$ and its truncation $G_r$ along a bandwidth of radius $r$ along the diagonal of the domain.

\subsection{Multipole graph neural operators} \label{sec_multipole}

Multipole graph neural operator (MGNO) has been introduced by~\citep{li2020multipole} and is a class of multi-scale networks that extends the graph neural operator architecture described in~\cref{sec_local} to capture the long-range interactions of the Green kernel. The main idea behind MGNO is to decompose the Green kernel $G$ into a sum of low-rank kernels as $G = K_1+\cdots+K_L$, which approximates the short and wide-range interactions in the PDEs. This architecture is motivated by the same reasons that led to the development of hierarchical low-rank matrices (see~\cref{sec_hierarchical}), such as the fast multipole method~\citep{greengard_rokhlin_1997,ying2004kernel}. It allows for the evaluation of the integral operation in \cref{eq_layer_update} in linear complexity.

MGNO is based on low-rank approximations of kernels, similar to DeepONets or low-rank neural operators (see~\cref{sec_deeponets} and~\citealt{li2020neural}), but is more flexible than vanilla DeepONets since it does not require the underlying kernels to be low-rank. Hence, if we consider a Green's function $G$ associated with a uniformly elliptic PDE in the form of \cref{eq_div_form}, then Weyl's law \citep{Weyl1911asymptotische,Canzani2013analysis, Minakshisundaram1949some} states that the eigenvalues of the solution operator associated with \cref{eq_div_form} decay at an algebraic rate of $\lambda_n \sim c n^{-2/d}$ for a constant $c>0$. This implies that the approximation error between the solution operator and its best rank-$k$ approximant decays only algebraically with $k$. Moreover, the decay rate deteriorates in high dimensions. In particular, the length $p$ of the feature vector in DeepONets must be significantly large to approximate the solution operator to a prescribed accuracy.

\begin{figure}[htbp]
    \centering
    \begin{overpic}[width=\textwidth]{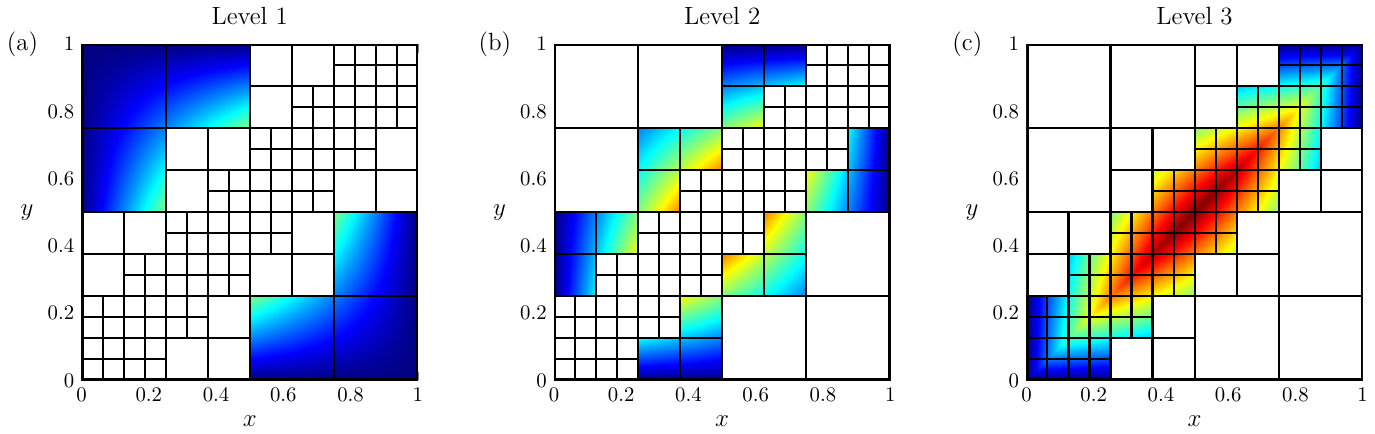}
    \end{overpic}
    \caption{Decomposition of the Green's function associated with the 1D Poisson equation, displayed in \cref{fig_green_diagonal}(a), into a hierarchy of kernels capturing different ranges of interactions: from long-range (a) to short-range (c) interactions.}
    \label{fig_green_hierarchy}
\end{figure}

However, \citet[Thm.~2.8]{bebendorf2003existence} showed that the Green's function $G$ associated with \cref{eq_div_form} can be well approximated by a low-rank kernel when restricted to separated subdomains $D_X\times D_Y$ of $\Omega\times \Omega$, satisfying the strong admissibility condition: $\dist(D_X,D_Y)<\diam(D_Y)$. Here, the distance and diameter in $\R^d$ are defined as
\[\dist(D_X,D_Y)=\inf_{x\in D_X, y\in D_Y}|x-y|,\qquad \diam(D_Y) = \sup_{y_1,y_2\in D_Y}|y_1-y_2|.\]
Then, for any $\epsilon\in(0,1)$, there exists a separable approximation of the form $G_k(x,y)=\sum_{i=1}^k u_i(x)v_i(y)$, with $k=\mathcal{O}(\log(1/\epsilon)^{d+1})$, such that
\[\|G-G_k\|_{L^2(D_X\times D_Y)}\leq \epsilon \|G\|_{L^2(D_X\times \hat{D}_Y)},\]
where $\hat{D}_Y$ is a domain slightly larger than $D_Y$~\citep[Thm.~2.8]{bebendorf2003existence}. This property has been exploited by~\citep{boulle2021learning,boulle2023elliptic,boulle2022learning} to derive sample complexity bounds for learning Green's functions associated with elliptic and parabolic PDEs. It motivates the decomposition of the Green kernel into a sum of low-rank kernels $G = K_1+\ldots+K_L$ in MGNO architectures. Indeed, one can exploit the low-rank structure of Green's functions on well-separated domains to perform a hierarchical decomposition of the domain $\Omega\times\Omega$ into a tree of subdomains satisfying the admissibility condition. In \cref{fig_green_hierarchy}, we illustrate the decomposition of the Green's function associated with the 1D Poisson equation on $\Omega=[0,1]$ with homogeneous Dirichlet boundary conditions into a hierarchy of $L=3$ levels of different range of interactions. The first level captures the long-range interactions, while the last level captures the short-range interactions. Then, the integral operation in \cref{eq_layer_update} can be performed by aggregating the contributions of the subdomains in the tree, starting from the leaves and moving up to the root. This allows for the evaluation of the integral operation in \cref{eq_layer_update} in linear complexity in the number of subdomains. The key advantage is that the approximation error on each subdomain decays exponentially fast as the rank of the approximating kernel increases.

One alternative approach to MGNO is to encode the different scales of the solution operators using a wavelet basis. This class of operator learning techniques~\citep{feliu2020meta,gupta2021multiwaveletbased,tripura2022wavelet} is based on the wavelet transform and aims to learn the solution operator kernel at multiple scale resolutions. One advantage over MGNO is that it does not require building a hierarchy of meshes, which could be computationally challenging in high dimensions or for complex domain geometries.

Finally, motivated by the success of the self-attention mechanism in transformers architectures for natural language processing~\citep{vaswani2017attention} and image recognition~\citep{dosovitskiy2020image}, several architectures have been proposed to learn global correlations in solution operators of PDEs. In particular, \citep{cao2021choose} introduced an architecture based on the self-attention mechanism for operator learning and observed higher performance on benchmark problems when compared against the Fourier Neural Operator. More recently, \citep{kissas2022learning} propose a Kernel-Coupled Attention mechanism to learn correlations between the entries of a vector feature representation of the output functions. In contrast, \citep{hao2023gnot} designed a general neural operator transformer (GNOT) that allows for multiple input functions and complex meshes.

\section{Learning neural operators} \label{sec_learning}

In this section, we discuss various technical aspects involved in training neural operators, such as the data acquisition of forcing terms and solutions, the amount of training data required in practice, and the optimization algorithms and loss functions used to train neural operators.

\subsection{Data acquisition} \label{sec_data_acquisition}

This section focuses on the data acquisition process for learning neural operators. In real-world applications, one may not have control over the distribution of source terms and solutions or locations of the sensors to measure the solutions at specific points in the domain. Therefore, we consider an idealized setting where one is interested in generating synthetic data using numerical PDE solvers to develop neural operator architectures. In this case, one has complete control over the distribution of source terms and solutions, as well as the locations of the sensors.

\subsubsection{Distribution of source terms} \label{sec_distribution}

The source terms $\{f_j\}_{j=1}^N$ used to generate pairs of training data to train neural operators are usually chosen to be random functions, sampled from a Gaussian random field~\citep{lu2021learning}. Let $\Omega\subset\R^d$ be a domain, then a stochastic process $\{X_x,\,x\in \Omega\}$ indexed by $\Omega$, is Gaussian if, for every finite set of indices $x_1,\ldots,x_n\in \Omega$, the vector of random variables $(X_{x_1},\ldots,X_{x_n})$ follows a multivariate Gaussian distribution. The Gaussian process (GP) distribution is completely determined by the following mean and covariance functions~\citep[Sec.~1.6]{adler2010geometry}:
\[\mu(x) = \mathbb{E}\{X_x\},\quad K(x,y) = \mathbb{E}\{[X_x-\mu(x)]^\top[X_y-\mu(y)]\},\quad x\in\Omega.\]
In the rest of the paper, we will denote a Gaussian process with mean $\mu$ and covariance kernel $K$ by $\mathcal{GP}(\mu,K)$. The mean function $\mu$ is usually chosen to be zero, while $K$ is symmetric and positive-definite.

When $K$ is continuous Mercer's theorem~\citep{mercer1909functions} states that there exists an orthonormal basis of eigenfunctions $\{\psi_j\}_{j=1}^\infty$ of $L^2(\Omega)$, and nonnegative eigenvalues $\lambda_1\geq \lambda_2\geq\cdots >0$ such that
\[
    K(x,y) = \sum_{j=1}^\infty\lambda_j\psi_j(x)\psi_j(y),\quad x,y\in \Omega,
\]
where the sum is absolutely and uniformly convergent~\citep[Thm.~4.6.5]{hsing2015theoretical}. Here, the eigenvalues and eigenfunctions of the kernel are defined as solutions to the associated Fredholm integral equation:
\[
    \int_\Omega K(x,y)\psi_j(y)\d y = \lambda_j\psi_j(x),\quad x\in \Omega.
\]
Then, the Karhunen--Lo\`eve theorem~\citep{karhunen1946lineare,loeve1946functions} ensures that a zero mean square-integrable Gaussian process $X_x$ with continuous covariance function $K$ admits the following representation:
\begin{equation} \label{eq_KL_expansion}
    X_x = \sum_{j=1}^\infty\sqrt{\lambda_j}c_j\psi_j(x),\quad c_j\sim\mathcal{N}(0,1),\quad x\in \Omega,
\end{equation}
where $c_j$ are independent and identically distributed (i.i.d.) random variables, and the convergence is uniform in $x\in \Omega$. Suppose the eigenvalue decomposition of the covariance function is known. In that case, one can sample a random function from the associated GP, $\mathcal{GP}(0,K)$, by sampling the coefficients $c_j$ in \cref{eq_KL_expansion} from a standard Gaussian distribution and truncated the series up to the desired resolution. Under suitable conditions, one can relate the covariance function $ K$'s smoothness to the random functions sampled from the associated GP~\citep[Sec.~3]{adler2010geometry}. Moreover, the decay rate of the eigenvalues provides information about the smoothness of the underlying kernel~\citep{ritter1995multivariate,zhu1997gaussian}. In practice, the number of eigenvalues greater than machine precision dictates the dimension of the finite-dimensional vector space spanned by the random functions sampled from $\mathcal{GP}(0,K)$.

\begin{figure}[htbp]
    \centering
    \begin{overpic}[width=\textwidth]{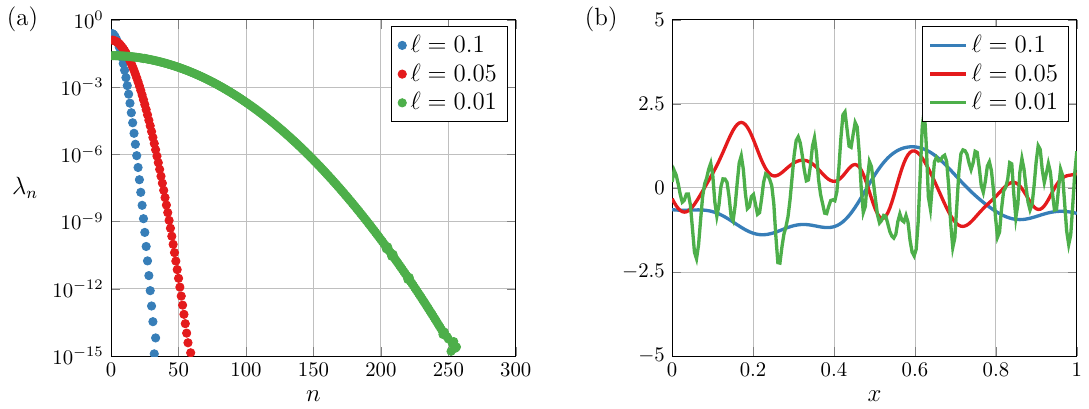}
    \end{overpic}
    \caption{(a) Eigenvalues of the squared-exponential kernel on $\Omega=[0,1]$ with length-scale parameters $\ell\in\{0.1,0.05,0.01\}$. (b) Random functions sampled from the associated Gaussian process $\mathcal{GP}(0,K)$, where $K$ is the squared-exponential kernel with length-scale parameters $\ell\in\{0.1,0.05,0.01\}$.}
    \label{fig_gp}
\end{figure}

One of the most common choices of covariance functions for neural operator learning include the squared-exponential kernel~\citep{lu2021learning,boulle2021data}, which is defined as
\[K(x,y) = \exp(-|x-y|^2/(2\ell^2)),\quad x,y\in \Omega,\]
where $\ell>0$ is the length-scale parameter, which roughly characterizes the distance at which two point values of a sampled random function become uncorrelated~\citep[Chapt.~5]{rasmussen2006gaussian}. Moreover, eigenvalues of the squared-exponential kernel decay exponentially fast at a rate that depends on the choice of $\ell$~\citep{zhu1997gaussian,boulle2022generalization}. After a random function $f$ has been sampled from the GP, one typically discretizes it by performing a piecewise linear interpolation at sensor points $x_1,\ldots,x_m\in \Omega$ by evaluating $f$ as these points. The interpolant can then solve the underlying PDE or train a neural operator. The number of sensors is chosen to resolve the underlying random functions and depends on their smoothness. Following the analysis by~\citet[Suppl.~Inf.~S4]{lu2021learning}, in one dimension, the error between $f$ and its piecewise linear interpolant is of order $\mathcal{O}(1/(m^2\ell^2))$, and one should choose $m\geq 1/\ell$. A typical value of $\ell$ lies in the range $\ell\in [0.01,0.1]$ with $m=100$ sensors~\citep{lu2021learning,boulle2021data}. We illustrate the eigenvalues of the squared-exponential kernel on $\Omega=[0,1]$ with length-scale parameters $\ell\in\{0.1,0.05,0.01\}$, along with the corresponding random functions sampled from the associated GP in \cref{fig_gp}. As the length-scale parameter $\ell$ decreases, the eigenvalues decay faster, and the sampled random functions become smoother.

Another possible choice of covariance functions for neural operator learning~\citep{benitez2023outofdistributional,zhu2023reliable} comes from the Mat\'ern class of covariance functions~\citep{rasmussen2006gaussian,stein1999interpolation}:
\[
    K(x,y) = \frac{2^{1-\nu}}{\Gamma(\nu)}\left(\frac{\sqrt{2\nu}|x-y|}{\ell}\right)^\nu K_\nu\left(\frac{\sqrt{2\nu}|x-y|}{\ell}\right),\quad x,y\in \Omega,
\]
where $\Gamma$ is the Gamma function, $K_\nu$ is a modified Bessel function, and $\nu$, $\ell$ are positive parameters that enable the control of the smoothness of the sampled random functions. Hence, the resulting Gaussian process is $\lceil\nu\rceil-1$ times mean-squared differentiable~\citep[Sec.~4.2]{rasmussen2006gaussian}. Moreover, the Mat\'ern kernel converges to the squared-exponential covariance function as $\nu\to\infty$. We refer to the book by \citet{rasmussen2006gaussian} for a detailed analysis of other standard covariance functions in Gaussian processes.

Finally, \citep{li2020fourier,kovachki2023neural} propose to use a Green kernel associated with a differential operator, which is a power of the Helmholtz equation, as a covariance function:
\[K = A(-\nabla^2+c I)^{-\nu}.\]
Here, $A$, $\nu>0$, and $c\geq 0$ are parameters that respectively govern the scaling of the Gaussian process, the algebraic decay rate of the spectrum, and the frequency of the first eigenfunctions of the covariance function. One motivation for this choice of distribution is that it allows the enforcement of prior information about the underlying model, such as the order of the differential operator, directly into the source terms. A similar behavior has been observed in a randomized linear algebra context when selecting the distribution of random vectors for approximating matrices from matrix-vector products using the randomized SVD~\citep{boulle2022generalization}. For example, the eigenvalues associated with this covariance kernel decay at an algebraic rate, implying that random functions sampled from this GP would be more oscillatory. This could lead to higher performance of the neural operator on high-frequency source terms. However, one downside of this approach is that a poor choice of the parameters can affect the training and approximation error of the neural operator. Hence, the studies~\citep{li2020fourier,kovachki2023neural} employ different parameter choices in each of the reported numerical experiments, suggesting that the covariance hyperparameters have been heavily optimized.

\subsubsection{Numerical PDE solvers} \label{sec_solvers}

After choosing the covariance function and generating the random source terms from the Gaussian process, one must solve the underlying PDE to generate the corresponding solutions. In general, the PDE is unknown, and one only has access to an oracle (such as physical experiments or black-box numerical solver) that outputs solutions $u$ to the PDE from input source terms $f$ as $\mathcal{L} u = f$. However, generating synthetic data from known mathematical models to design or evaluate neural operator architectures is often convenient. One can use numerical PDE solvers to generate the corresponding solutions in this case. This section briefly describes the different numerical methods that can be used to solve the underlying PDEs, along with their key attributes summarized in \cref{tab_solvers}. We want to emphasize that these methods have many variations, and we refer to the most standard ones.

\begin{table}[htbp]
    \centering
    \begin{tabular}{ p{4cm} c c c  }
        \hline
        Property         & Finite differences & Finite elements & Spectral methods \\
        \hline
        Domain geometry  & Simple             & Complex         & Simple           \\
        Approximation    & Local              & Local           & Global           \\
        Linear system    & Large sparse       & Large sparse    & Small dense      \\
        Convergence rate & Algebraic          & Algebraic       & Spectral         \\
        \hline
    \end{tabular}
    \caption{Summary of the different properties of standard finite difference, finite element, and spectral methods for solving PDEs.}
    \label{tab_solvers}
\end{table}

When the PDE does not depend on time, the most common techniques for discretizing and solving it are finite difference methods (FDM), finite element methods (FEM), and spectral methods. The finite difference method consists of discretizing the computational domain $\Omega$ into a grid and approximating spatial derivatives of the solution $u$ from linear combinations of the values of $u$ at the grid points using finite difference operators~\citep[Chap.~8]{iserles2009first}. This approach is based on a local Taylor expansion of the solution and is very easy to implement on rectangular geometries. However, it usually requires a uniform grid approximation of the domain, which might not be appropriate for complex geometries and boundary conditions. The finite element method~\citep[Chap.~14]{suli2003introduction} employs a different approach than FDM and considers the approximation of the solution $u$ on a finite-dimensional vector space spanned by basis functions with local support on $\Omega$. The spatial discretization of the domain $\Omega$ is performed via a mesh representation. The basis functions are often chosen as piecewise polynomials supported on a set of elements, which are adjacent cells in the mesh. This approach is more flexible than FDM and can be used to solve PDEs on complex geometries and boundary conditions. However, it is more challenging to implement and requires the construction of a mesh of the domain $\Omega$, which can be computationally expensive. We highlight that the finite difference and finite element methods lead to large, sparse, and highly structured linear algebra systems, which can be solved efficiently using iterative methods. Two commonly used finite element software for generating training data for neural operators are FEniCS~\citep{alnaes2015fenics} and Firedrake~\citep{rathgeber2016firedrake,ham2023firedrake}. These open-source software exploit the Unified Form Language (UFL) developed by~\citep{alnaes2014unified} to define the weak form of the PDE in a similar manner as in mathematics and automatically generate the corresponding finite element assembly code before exploiting fast linear and nonlinear solvers through a Python interface with the high-performance PETSc library~\citep{balay2019petsc}.

Finally, spectral methods~\citep{iserles2009first,gottlieb1977numerical,trefethen2000spectral} are based on the approximation of the solution $u$ on a finite-dimensional vector space spanned by basis functions with global support on $\Omega$, usually Chebyshev polynomials or trigonometric functions. Spectral methods are motivated by the fact that the solution $u$ of a PDE defined on a 1D interval is often smooth if the source term is smooth and so can be well-approximated by a Fourier series if it is periodic or a Chebyshev series if it is not periodic. Hence, spectral methods lead to exponential convergence, also called spectral accuracy, to analytic solutions with respect to the number of basis functions, unlike FDM and FEM, which only converge at an algebraic rate. The fast convergence rate of spectral methods implies that a small number of basis functions is usually required to achieve a given accuracy. Therefore, the matrices associated with the resulting linear algebra systems are much smaller than for FEM but are dense. In summary, spectral methods have a competitive advantage over FDM and FEM on simple geometries, such as tensor-product domains, and when the solution is smooth. At the same time, FEM might be difficult to implement but is more flexible. A convenient software for solving simple PDEs using spectral methods is the Chebfun software system~\citep{driscoll2014chebfun}, an open-source package written in MATLAB.

For time-dependent PDEs, one typically starts by performing a time-discretization using a time-stepping scheme, such as backward differentiation schemes (e.g.~backward Euler) and Runge--Kutta methods~\citep{iserles2009first,suli2003introduction}, and then employ a spatial discretization method, such as the techniques described before in this section, to solve the resulting stationary PDE at each time-step.

\subsubsection{Amount of training data} \label{sec_number}

Current neural operator approaches typically require a relatively small amount of training data, in the order of a thousand input-output pairs, to approximate solution operators associated with PDEs~\citep{lu2021learning,goswami2023physics,kovachki2023neural,boulle2023elliptic}. This contrasts with the vast amount of data used to train neural networks for standard supervised learning tasks, such as image classification, which could require hundreds of millions of labeled samples~\citep{lecun2015deep}. This difference can be explained by the fact that solution operators are often highly structured, which can be exploited to design data-efficient neural operator architectures (see~\cref{sec_architectures}).

Recent numerical experiments have shown that the rate of convergence of neural operators with respect to the number of training samples evolves in two regimes~\citep[Fig.~S10]{lu2021learning}. In the first regime, we observe a fast decay of the testing error at an exponential rate~\citep{boulle2023elliptic}. Then, the testing error decays at a slower algebraic rate in the second regime for a larger amount of samples and saturates due to discretization error and optimization issues, such as convergence to a suboptimal local minima.

On the theoretical side, several works derived sample complexity bounds that characterize the amount of training data required to learn solution operators associated with certain classes of linear PDEs to within a prescribed accuracy $0<\epsilon<1$. In particular, \citep{boulle2023elliptic,schafer2021sparse} focus on approximating Green's functions associated with uniformly elliptic PDEs in divergence form:
\begin{equation} \label{eq_div_form_2}
    -\vdiv(A(x)\nabla u) = f,\quad x\in \Omega\subset \R^d,
\end{equation}
where $A(x)$ is a symmetric bounded coefficient matrix (see~\cref{eq_div_form}). These studies construct data-efficient algorithms that converge exponentially fast with respect to the number of training pairs. Hence, they can recover the Green's function associated with \cref{eq_div_form_2} to within $\epsilon$ using only $\mathcal{O}(\text{polylog}(1/\epsilon))$ sample pairs. The method employed by \citep{boulle2023elliptic} consists of recovering the hierarchical low-rank structure satisfied by Green's function on well-separated subdomains~\citep{bebendorf2003existence,lin2011fast,levitt2022randomized} using a generalization of the rSVD to Hilbert--Schmidt operators~\citep{boulle2022generalization,boulle2021learning,halko2011finding,martinsson2020randomized}. Interestingly, the approach by \citep{schafer2021sparse} is not based on low-rank techniques but relies on the sparse Cholesky factorization of elliptic solution operators~\citep{schafer2021compression}.

Some of the low-rank recovery techniques employed by~\citep{boulle2023elliptic} extend naturally to time-dependent parabolic PDEs~\citep{boulle2022learning} in spatial dimension $d\geq 1$ of the form:
\begin{equation} \label{eq_parabolic}
    \frac{\partial u}{\partial t} - \vdiv(A(x,t)\nabla u) = f,\quad x\in \Omega\subset \R^d,\quad t\in (0,T],
\end{equation}
where the coefficient matrix $A(x,t)\in \R^{d\times d}$ is symmetric positive definite with bounded coefficient functions in $L^\infty(\Omega\times [0,T])$, for some $0<T<\infty$, and satisfies the uniform parabolicity condition~\citep[Sec.~7.1.1]{evans2010partial}. Parabolic systems in the form of \cref{eq_parabolic} model various time-dependent phenomena, including heat conduction and particle diffusion. \citet[Thm.~9]{boulle2022learning} showed that the Green's function associated with \cref{eq_parabolic} admits a hierarchical low-rank structure on well-separated subdomains, similarly to the elliptic case~\citep{bebendorf2003existence}. Combining this with the pointwise bounds satisfied by the Green's function~\citep{cho2012global}, one can construct an algorithm that recovers the Green's function to within $\epsilon$ using $\mathcal{O}(\text{poly}(1/\epsilon))$ sample pairs~\citep[Thm.~10]{boulle2022learning}.

Finally, other approaches~\citep{de2021convergence,jin2022minimax,stepaniants2023learning} derived convergence rates for a broader class of operators between infinite-dimensional Hilbert spaces, which are not necessarily associated with solution operators of PDEs. In particular, \citep{de2021convergence} consider the problem of estimating the eigenvalues of an unknown, and possibly unbounded, self-adjoint operator assuming that the operator is diagonalizable in a known basis of eigenfunctions, and highlight the impact of varying the smoothness of training and test data on the convergence rates. Then, \citep{stepaniants2023learning,jin2022minimax} derive upper and lower bounds on the sample complexity of Hilbert--Schmidt operators between two reproducing kernel Hilbert spaces (RKHS) that depend on the smoothness of the input and output functions.

\subsection{Optimization} \label{sec_optimization}

Once the neural operator architecture has been selected and the training dataset is constituted, the next task is to train the neural operator by solving an optimization problem in the form of \cref{eq_loss_function}. The aim is to identify the best parameters of the underlying neural network so that the output $\hat{\mathcal{A}}(f;\theta)$ of the neural operator evaluated at a forcing term $f$ in the training dataset fits the corresponding ground truth solution $u$. This section describes the most common choices of loss functions and optimization algorithms employed in current operator learning approaches. Later in \cref{sec_convergence}, we briefly discuss how to measure the convergence and performance of a trained neural operator.

\subsubsection{Loss functions} \label{sec_loss}

The choice of loss function in operator learning is a critical step, as it directs the optimization process and ultimately affects the model's performance. Different types of loss functions can be utilized depending on the task's nature, the operator's structure, and the function space's properties. A common choice of loss function in ML is the mean squared error (MSE), which is defined as
\begin{equation} \label{eq_mse}
    L_{\text{MSE}} = \frac{1}{N}\sum_{i=1}^N\frac{1}{m}\sum_{j=1}^m |\hat{A}(f_i)(x_j)-u_i(x_j)|^2\approx \frac{1}{N}\sum_{i=1}^N\|\hat{A}(f_i)-u_i\|_{L^2(\Omega)}^2,
\end{equation}
and is employed in the original DeepONet study~\citep{lu2021learning}. Here, $N$ is the number of training samples, $m$ is the number of sensors, $f_i$ is the $i$-th forcing term, $u_i$ is the corresponding ground truth solution, $\hat{A}(f_i)$ is the output of the neural operator evaluated at $f_i$, and $x_j$ is the $j$-th sensor location. This loss function discretizes the squared $L^2$ error between the output of the neural operator and the ground truth solution at the sensor locations. When the sensor grid is regular, one can employ a higher-order quadrature rule to discretize the $L^2$ norm. Moreover, in most cases, it may be beneficial to use a relative error, especially when the magnitudes of the output function can vary widely. Then, \citep{boulle2021data} use the following relative squared $L^2$ loss function
\begin{equation} \label{eq_relative_l2}
    L = \frac{1}{N}\sum_{i=1}^N\frac{\|\hat{A}(f_i)-u_i\|_{L^2(\Omega)}^2}{\|u_i\|_{L^2(\Omega)}^2},
\end{equation}
which is discretized using a trapezoidal rule. The most common loss function in operator learning is the relative $L^2$ error employed in Fourier neural operator techniques~\citep{li2020fourier}:
\begin{equation} \label{eq_l2}
    L_2 = \frac{1}{N}\sum_{i=1}^N\frac{\|\hat{A}(f_i)-u_i\|_{L^2(\Omega)}}{\|u_i\|_{L^2(\Omega)}}.
\end{equation}
\citep{kovachki2023neural} observed a better normalization of the model when using a relative loss function, and that the choice of \cref{eq_l2} decreases the testing error by a factor of two compared to \cref{eq_relative_l2}.

For tasks that require robustness to outliers or when it is important to measure the absolute deviation, the $L^1$ loss can be employed. It is defined as
\begin{equation} \label{eq_l1}
    L_1 = \frac{1}{N}\sum_{i=1}^N\frac{\|\hat{A}(f_i)-u_i\|_{L^1(\Omega)}}{\|u_i\|_{L^1(\Omega)}}.
\end{equation}
This loss function tends to be less sensitive to large deviations than the $L^2$ loss~\citep{alpak2023augmenting,lyu2023multi,zhao2024recfno}. Furthermore, Sobolev norms can also be used as a loss function when the unknown operator $\mathcal{A}$ involves functions in Sobolev spaces~\citep[Chapt.~5]{evans2010partial}, particularly when the derivatives of the input and output functions play a role~\citep{son2021sobolev,yu2023tuning,o2024derivative}. For example, one could perform training with a relative $H^1$ loss to enforce the smoothness of the neural operator output. Finally, when the underlying PDE is known, one can enforce it as a weak constraint when training the neural operator by adding a PDE residual term to the loss function~\citep{li2021physics,wang2021learning}, similarly to physics-informed neural networks~\citep{raissi2019physics}.

\subsubsection{Optimization algorithms and implementation} \label{sec_algorithms}

The training procedure of neural operators is typically performed using Adam optimization algorithm~\citep{kingma2015adam,kovachki2023neural,lu2021learning,li2020fourier,goswami2023physics} or one of its variants such as AdamW~\citep{loshchilov2019decoupled,hao2023gnot}. Hence, the work introducing DeepONets by \citep{lu2021learning} employed Adam algorithm to train the neural network architecture with a default learning rate of $0.001$. In contrast, \citep{kovachki2023neural} incorporate learning rate decay throughout the optimization of Fourier neural operators. One can also employ a two-step training approach by minimizing the loss function using Adam algorithm for a fixed number of iterations and then fine-tuning the neural operator using the L-BFGS algorithm~\citep{byrd1995limited,cuomo2022scientific,he2020physics,mao2020physics,boulle2021data}. This approach has been shown to improve the convergence rate of the optimization when little data is available in PINNs applications~\citep{he2020physics}. Popular libraries for implementing and training neural operators include PyTorch~\citep{paszke2019pytorch} and TensorFlow~\citep{abadi2016tensorflow}.

Thus far, there has been limited focus on the theoretical understanding of convergence and optimization of neural operators. Since neural operators are a generalization of neural networks in infinite dimensions, existing convergence results of physics-informed neural networks~\citep{wang2021eigenvector,wang2022and} based on the neural tangent kernel (NTK) framework~\citep{jacot2018neural,du2019gradient,allen2019convergence} should naturally extend to neural operators. One notable exception is the study by \citep{wang2022improved}, which analyzes the training of physics-informed DeepONets~\citep{wang2021learning} and derives a weighting scheme guided by NTK theory to balance the data and the PDE residual terms in the loss function.

\subsubsection{Measuring convergence and super-resolution} \label{sec_convergence}

After training a neural operator, one typically measures its performance by evaluating the testing error, such as the relative $L^2$-error, on a set of unseen data generated using the procedure described in \cref{sec_learning}. In general, state-of-the-art neural operator architectures report a relative testing error of $1\%-10\%$ depending on the problems considered~\citep{kovachki2023neural,lu2021learning,li2020fourier}.

However, it is essential to note that the testing error may not be a good measure of the performance of a neural operator, as it does not provide any information about the generalization properties of the model. Hence, the testing forcing terms are usually sampled from the same distribution as the training forcing terms, so they lie on the same finite-dimensional function space, determined by the spectral decay of the GP covariance kernel eigenvalues (see~\cref{sec_data_acquisition}). Moreover, in real applications, the testing source terms could have different distributions than the ones used for training ones, and extrapolation of the neural operator may required. \citep{zhu2023reliable} investigate the extrapolation properties of DeepONet with respect to the length-scale parameter of the underlying source term GP. They observe that the testing error increases when the length-scale parameter corresponding to the test data decreases. At the same time, the neural operator can extrapolate to unseen data with a larger length scale than the training dataset, i.e., smoother functions.

\begin{figure}[htbp]
    \centering
    \begin{overpic}[width=\textwidth]{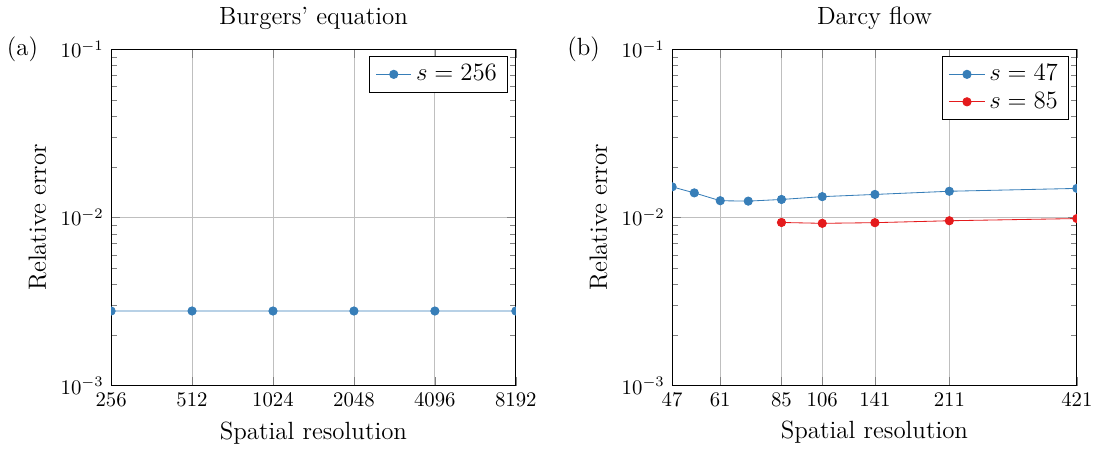}
    \end{overpic}
    \caption{(a) Relative test error at different spatial resolutions of the Fourier neural operator trained to approximate the solution operator of the 1D Burgers' equation~\eqref{eq_burgers} with trained resolution of $s=256$. (b) FNO trained on 2D Darcy flow a resolution of $s=47$ (blue) and $s=85$ (red) and evaluated at higher spatial resolutions.}
    \label{fig_super_resolution}
\end{figure}

One attractive property of neural operators is their resolution invariance to perform predictions at finer spatial resolutions than the training dataset on which they have been trained. This is usually called zero-shot super-resolution~\citep[Sec.~7.2.3]{kovachki2023neural}. To investigate this property, we reproduce the numerical examples of~\citep[Sec.~7.2]{kovachki2023neural} and train a Fourier neural operator to approximate the solution operator of Burgers' equation and Darcy flow at a low-resolution data and evaluate the operator at higher resolutions. We consider the one-dimensional Burgers' equation:
\begin{equation} \label{eq_burgers}
    \frac{\partial }{\partial t}u(x,t) + \frac{1}{2}\frac{\partial}{\partial_x}(u(x,t)^2) = \nu\frac{\partial^2}{\partial_{x^2}}u(x,t),\quad x\in (0,2\pi),\quad t\in [0,1],
\end{equation}
with periodic boundary conditions and viscosity $\nu=0.1$. We are interested in learning the solution operator $\mathcal{A}:L^2_{\text{per}}((0,2\pi))\to H_{\text{per}}^1((0,2\pi))$, which maps initial conditions $u_0\in L^2_{\text{per}}((0,2\pi))$ to corresponding solutions $u(\cdot,1)\in H_{\text{per}}^1((0,2\pi))$ to \cref{eq_burgers} at time $t=1$. We then discretize the source and solution training data on a uniform grid with spatial resolution $s=256$ and evaluate the trained neural operator at finer spatial resolutions $s\in\{512,1024,2048\}$. We observe in \cref{fig_super_resolution}(a) that the relative testing error of the neural operator is independent of the spatial resolutions, as reported by~\citet[Sec.~7.2]{kovachki2023neural}.

\begin{figure}[htbp]
    \centering
    \begin{overpic}[width=\textwidth]{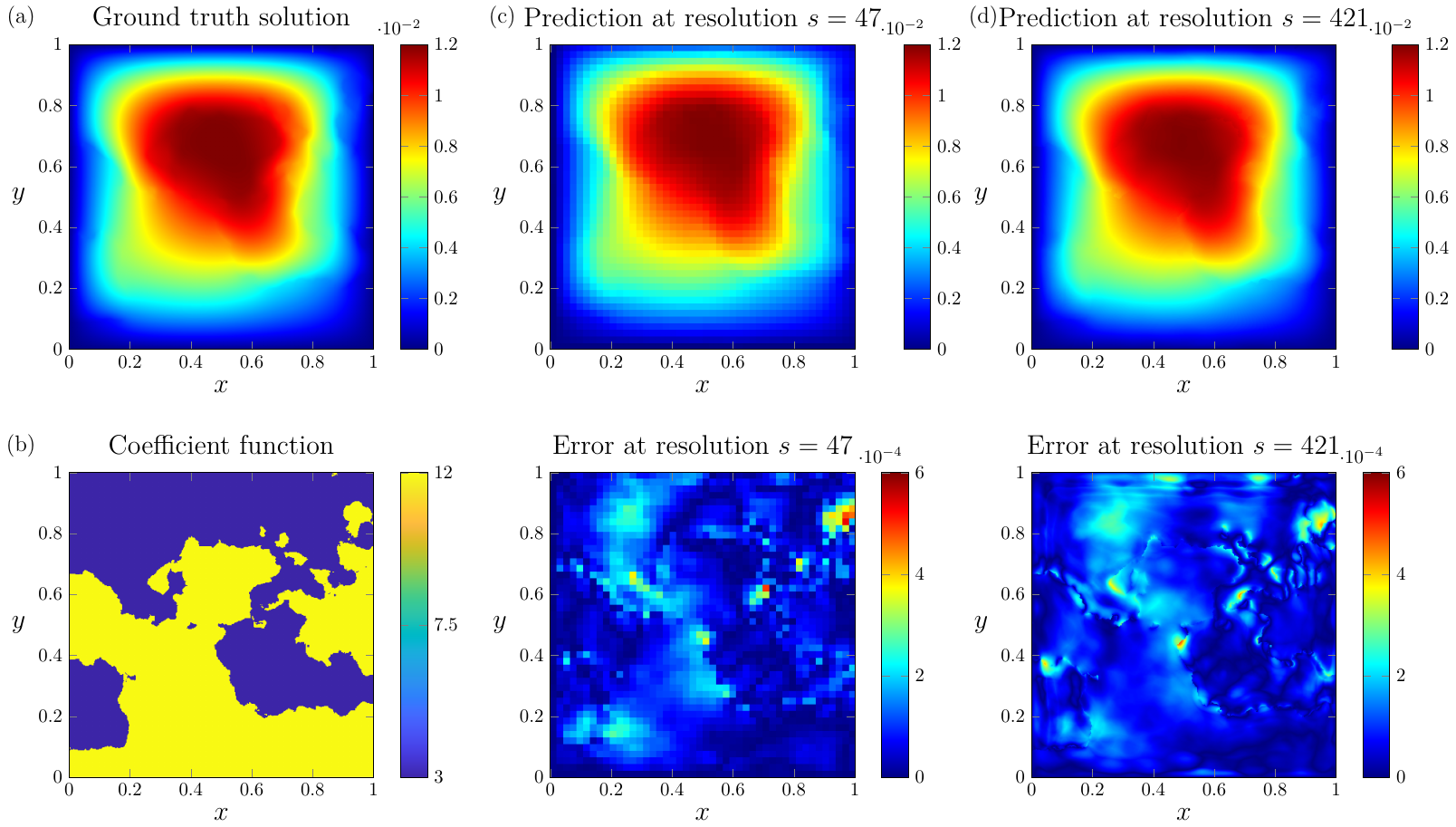}
    \end{overpic}
    \caption{(a) Ground truth solution $u$ to the 2D Darcy flow equation~\eqref{eq_darcy} corresponding to the coefficients function plotted in (b). (c)-(d) Predicted solution and approximation error at $s=47$ and $s=421$ by a Fourier neural operator trained on a Darcy flow dataset with spatial resolution of $s=47$.}
    \label{fig_super_resolution_image}
\end{figure}

Next, we consider the two-dimensional Darcy flow equation~\eqref{eq_darcy} with constant source term $f=1$ and homogeneous Dirichlet boundary conditions on a unit square domain $\Omega = [0,1]^2$:
\begin{equation} \label{eq_darcy}
    -\vdiv(a(x)\nabla u) = 1,\quad x\in [0,1]^2.
\end{equation}
We train a Fourier neural operator to approximate the solution operator, mapping the coefficient function $a$ to the associated solution $u$ to \cref{eq_darcy}. We reproduce the numerical experiment in \citep[Sec.~6.2]{kovachki2023neural}, where the random coefficient functions $a$ are piecewise constant. The random functions $a$ are generated as $a\sim T\circ f$, where $f\sim \mathcal{GP}(0,C)$, with $C = (-\Delta+9I)^{-2}$ and $T:\R\to\R^{+}$ is defined as
\[
    T(x) =
    \begin{cases}
        12, & \text{if } x\geq 0, \\
        3,  & \text{if } x<0.
    \end{cases}
\]
We discretize the coefficient and solution training data on a $s\times s$ uniform grid with spatial resolution $s=47$ and evaluate the trained neural operator at higher spatial resolutions in \cref{fig_super_resolution}. The relative testing error does not increase as the spatial resolution increases. Moreover, training the neural operator on a higher spatial resolution dataset can decrease the testing error. We also plot the ground truth solution $u$ to \cref{eq_darcy} in \cref{fig_super_resolution_image}(a) corresponding to the coefficient function plotted in panel (b), along with the predicted solutions and approximation errors at $s=47$ and $s=421$ by the Fourier neural operator in panels (c) and (d). We want to point the reader interested in the discretization properties of neural operators to the recent perspective on representation equivalent neural operators (ReNO) by~\citep{bartolucci2023neural}.

\section{Conclusions and future challenges} \label{sec_physics}

In this paper, we provided a comprehensive overview of the recent developments in neural operator learning, a new paradigm at the intersection of scientific computing and ML for learning solution operators of PDEs. Given the recent surge of interest in this field, a key question concerns the choice of neural architectures for different PDEs. Most theoretical studies in the field analyze and compare neural operators through the prism of approximation theory. We proposed a framework based on numerical linear algebra and matrix recovery problems for interpreting the type of neural operator architectures that can be used to learn solution operators of PDEs. Hence, solution operators associated with linear PDEs can often be written as integral operators with a Green's function and recover by a one-layer neural operator, which after discretization is equivalent to a matrix recovery problem.

Moreover, the choice of architectures, such as FNO or DeepONet, enforces or preserves different properties of the PDE solution operator, such as being translation invariant, low-rank, or off-diagonal low-rank (see~\cref{tab_property_kernels}). We then focused on the data acquisition process. We highlighted the importance of the distribution of source terms, usually sampled from a Gaussian process with a tailored covariance kernel, on the resulting performance of the neural operator. Following recent works on elliptic and parabolic PDEs and numerical experiments, we also discussed the relatively small amount of training data needed for operator learning. Finally, we studied the different choices of optimization algorithms and loss functions and highlighted the super-resolution properties of neural operators, i.e., their ability to be evaluated at higher resolution than the training dataset with a minor impact on the performance. There are, however, several remaining challenges in the field.

\paragraph{Distribution of probes.} Most applications of neural operators employ source terms that are globally supported on the domain, sampled from a Gaussian process, and whose distribution is fixed before training. However, this might now always apply to real-world engineering or biological systems, where source terms could be localized in space and time. A significant problem is to study the impact of the distribution of locally supported source terms on the performance of neural operators, both from a practical and theoretical viewpoint. Hence, recent sample complexity works on elliptic and parabolic PDEs exploit structured source terms~\citep{boulle2023elliptic,boulle2022learning,schafer2021sparse}. Another area of future research is to employ adaptive source terms to fine-tune neural operators for specific applications. This could lead to higher performance by selecting source terms that maximize the training error or allow efficient transfer learning between different applications without retraining a large neural operator.

\paragraph{Software and datasets.}
An essential step towards democratizing operator learning involves the development of open-source software and datasets for training and comparing neural operators, similar to the role played by the MNIST~\citep{lecun1998gradient} and ImageNet~\citep{deng2009imagenet} databases in the improvement of computer vision techniques. However, due to the fast emerging methods in operator learning, there have been limited attempts beyond~\citep{lu2022comprehensive} to standardize the datasets and software used in the field. Establishing a list of standard PDE problems across different scientific fields, such as fluid dynamics, quantum mechanics, and epidemiology, with other properties (e.g.~linear/nonlinear, steady/time-dependent, low/high dimensional, smooth/rough solutions, simple/complex geometry) would allow researchers to compare and identify the neural operator architectures that are the most appropriate for a particular task. A recent benchmark has been proposed to evaluate the performance of physics-informed neural networks for solving PDEs~\citep{hao2023pinnacle}.

\paragraph{Real-world applications.} Neural operators have been successfully applied to perform weather forecasting and achieve spectacular performance in terms of accuracy and computational time to solutions compared to traditional numerical weather prediction techniques while being trained on historical weather data~\citep{kurth2023fourcastnet,lam2023learning}. An exciting development in the field of operator learning would be to expand the scope of applications to other scientific fields and train the models on real datasets, where the underlying PDE governing the data is unknown to discover new physics.

\paragraph{Theoretical understanding.}
Following the recent works on the approximation theory of neural operators and sample complexity bounds for different classes of PDEs, there is a growing need for a theoretical understanding of convergence and optimization. In particular, an exciting area of research would be to extend the convergence results of physics-informed neural networks and the neural tangent kernel framework to neural operators. This would enable the derivation of rigorous convergence rates for different types of neural operator architectures and loss functions and new schemes for initializing the weight distributions in the underlying neural networks.

\paragraph{Physical properties.}
Most neural operator architectures are motivated by obtaining a good approximation of the solution operator of a PDE. However, the resulting neural operator is often highly nonlinear, difficult to interpret mathematically, and might not satisfy the physical properties of the underlying PDE, such as conservation laws or symmetries~\citep{olver2000applications}. There are several promising research directions in operator learning related to symmetries and conservation laws~\citep{otto2023unified}. One approach would be to enforce known physical properties when training neural operators, either strongly through structure preserving architectures~\citep{richter2022neural}, or weakly by adding a residual term in the loss function~\citep{li2021physics,wang2021learning}. Another direction is to discover new physical properties of the underlying PDEs from the trained neural operator. While \citep{boulle2021data} showed that symmetries of linear PDEs can be recovered from the learned Green's function, this approach has not been extended to nonlinear PDEs. Finally, one could also consider using reinforcement learning techniques for enforcing physical constraints after the optimization procedure, similar to recent applications in large language models~\citep{ouyang2022training}.

\acks{The work of both authors was supported by the Office of Naval Research (ONR), under grant N00014-23-1-2729. N.B. was supported by an INI-Simons Postdoctoral Research Fellowship. A.T. was supported by National Science Foundation grants DMS-2045646 and a Weiss Junior Fellowship Award.}

\vskip 0.2in
\addcontentsline{toc}{section}{References}
\bibliography{biblio}

\end{document}